\documentclass[12pt]{article}
\usepackage{latexsym, amssymb}

\DeclareMathAlphabet\gothic{U}{euf}{m}{n}

\makeatletter
\def\eqnarray{\stepcounter{equation}\let\@currentlabel=\theequation
\global\@eqnswtrue
\tabskip\@centering\let\\=\@eqncr
$$\halign to \displaywidth\bgroup\hfil\global\@eqcnt\z@
  $\displaystyle\tabskip\z@{##}$&\global\@eqcnt\@ne
  \hfil$\displaystyle{{}##{}}$\hfil
  &\global\@eqcnt\tw@ $\displaystyle{##}$\hfil
  \tabskip\@centering&\llap{##}\tabskip\z@\cr}

\def\endeqnarray{\@@eqncr\egroup
      \global\advance\c@equation\m@ne$$\global\@ignoretrue}
\makeatother

\begin{document}
\bibliographystyle{tom}

\newtheorem{lemma}{Lemma}[section]
\newtheorem{thm}[lemma]{Theorem}
\newtheorem{cor}[lemma]{Corollary}
\newtheorem{voorb}[lemma]{Example}
\newtheorem{rem}[lemma]{Remark}
\newtheorem{prop}[lemma]{Proposition}
\newtheorem{stat}[lemma]{{\hspace{-5pt}}}
\newtheorem{obs}[lemma]{Observation}

\newenvironment{remarkn}{\begin{rem} \rm}{\end{rem}}
\newenvironment{exam}{\begin{voorb} \rm}{\end{voorb}}
\newtheorem{rems}[lemma]{Remarks}
\newenvironment{remarkns}{\begin{rems} \rm}{\end{rems}}

\newenvironment{obsn}{\begin{obs} \rm}{\end{obs}}
\newenvironment{obsns}{\begin{obss} \rm}{\end{obss}}

\newtheorem{obss}[lemma]{Observations}

\newcommand{\gota}{\gothic{a}}
\newcommand{\gotb}{\gothic{b}}
\newcommand{\gotc}{\gothic{c}}
\newcommand{\gote}{\gothic{e}}
\newcommand{\gotf}{\gothic{f}}
\newcommand{\gotg}{\gothic{g}}
\newcommand{\gothh}{\gothic{h}}
\newcommand{\gotk}{\gothic{k}}
\newcommand{\gotm}{\gothic{m}}
\newcommand{\gotn}{\gothic{n}}
\newcommand{\gotp}{\gothic{p}}
\newcommand{\gotq}{\gothic{q}}
\newcommand{\gotr}{\gothic{r}}
\newcommand{\gots}{\gothic{s}}
\newcommand{\gotu}{\gothic{u}}
\newcommand{\gotv}{\gothic{v}}
\newcommand{\gotw}{\gothic{w}}
\newcommand{\gotz}{\gothic{z}}
\newcommand{\gotG}{\gothic{G}}
\newcommand{\gotL}{\gothic{L}}
\newcommand{\gotS}{\gothic{S}}
\newcommand{\gotT}{\gothic{T}}

\newcounter{teller}
\renewcommand{\theteller}{\Roman{teller}}
\newenvironment{tabel}{\begin{list}%
{\rm \bf \Roman{teller}.\hfill}{\usecounter{teller} \leftmargin=1.1cm
\labelwidth=1.1cm \labelsep=0cm \parsep=0cm}
                      }{\end{list}}

\newcounter{tellerr}
\renewcommand{\thetellerr}{\roman{tellerr}}
\newenvironment{subtabel}{\begin{list}%
{\rm  \roman{tellerr}.\hfill}{\usecounter{tellerr} \leftmargin=1.1cm
\labelwidth=1.1cm \labelsep=0cm \parsep=0cm}
                         }{\end{list}}

\newcounter{proofstep}
\newcommand{\nextstep}{\refstepcounter{proofstep}\ruimte \par 
          \noindent{\bf Step \theproofstep} \hspace{5pt}}
\newcommand{\firststep}{\setcounter{proofstep}{0}\nextstep}

\newcommand{\Ni}{{\bf N}}
\newcommand{\Ri}{{\bf R}}
\newcommand{\Ci}{{\bf C}}
\newcommand{\Ti}{{\bf T}}
\newcommand{\Zi}{{\bf Z}}
\newcommand{\Fi}{{\bf F}}

\newcommand{\proof}{\mbox{\bf Proof} \hspace{5pt}} 
\newcommand{\remark}{\mbox{\bf Remark} \hspace{5pt}}
\newcommand{\ruimte}{\vskip10.0pt plus 4.0pt minus 6.0pt}

\newcommand{\ad}{{\mathop{\rm ad}}}
\newcommand{\Ad}{{\mathop{\rm Ad}}}
\newcommand{\Aut}{\mathop{\rm Aut}}
\newcommand{\arccot}{\mathop{\rm arccot}}
\newcommand{\clos}{{\mathop{\rm cl}}}

\newcommand{\diam}{\mathop{\rm diam}}
\newcommand{\divv}{\mathop{\rm div}}
\newcommand{\codim}{\mathop{\rm codim}}
\newcommand{\interior}{\mathop{\rm int}}
\newcommand{\RRe}{\mathop{\rm Re}}
\newcommand{\IIm}{\mathop{\rm Im}}
\newcommand{\Tr}{{\mathop{\rm Tr}}}
\newcommand{\supp}{\mathop{\rm supp}}
\newcommand{\sgn}{\mathop{\rm sgn}}
\newcommand{\esssup}{\mathop{\rm ess\,sup}}
\newcommand{\Int}{\mathop{\rm Int}}
\newcommand{\Leibniz}{\mathop{\rm Leibniz}}
\newcommand{\lcm}{\mathop{\rm lcm}}
\newcommand{\mod}{\mathop{\rm mod}}
\newcommand{\spann}{\mathop{\rm span}}
\newcommand{\ubar}{\underline{\;}}
\newcommand{\one}{1\hspace{-4.5pt}1}

\hyphenation{groups}
\hyphenation{unitary}

\newcommand{\dm}{{d_{\mbox{\gothicss m}}}}
\newcommand{\dms}{{d_{\mbox{\gothics m}}}}
\newcommand{\hc}{\overline{H}}
\newcommand{\sodp}{ \{ 1,\ldots,d' \} }
\newcommand{\ccig}{C_c^\infty(G)}
\newcommand{\Cprime}[1]{C^{#1}{\hspace{0.8pt}}^\prime}
\newcommand{\CCprime}[1]{{\cal C}^{#1}{\hspace{0.8pt}}^\prime}
\newcommand{\Cdeltaprime}[1]{C_\Delta^{#1}{\hspace{0.5pt}}^\prime}

\newcommand{\Gdual}{\widetilde{G}}
\newcommand{\htozp}{H'_{2;1}\raisebox{10pt}[0pt][0pt]{\makebox[0pt]{\hspace{-28pt}$\scriptstyle\circ$}}}
\newcommand{\htozpu}{H^{\prime(u)}_{2;1}\raisebox{10pt}[0pt][0pt]{\makebox[0pt]{\hspace{-36pt}$\scriptstyle\circ$}}}

\newcommand{\htozpud}{H^{(u)}_{2;1}\raisebox{10pt}[0pt][0pt]{\makebox[0pt]{\hspace{-32pt}$\scriptstyle\circ$}}}
\newcommand{\htozpukd}{H^{(u_k)}_{2;1}\raisebox{10pt}[0pt][0pt]{\makebox[0pt]{\hspace{-41pt}$\scriptstyle\circ$}}}
\newcommand{\htozpund}{H^{(u_n)}_{2;1}\raisebox{10pt}[0pt][0pt]{\makebox[0pt]{\hspace{-41pt}$\scriptstyle\circ$}}}

\newcommand{\htozpH}{H^{\prime(H)}_{2;1}\raisebox{10pt}[0pt][0pt]{\makebox[0pt]{\hspace{-36pt}$\scriptstyle\circ$}}}
\newcommand{\htozpun}{H^{\prime(u_n)}_{2;1}\raisebox{10pt}[0pt][0pt]{\makebox[0pt]{\hspace{-45pt}$\scriptstyle\circ$}}}
\newcommand{\htoz}{H_{2;1}\raisebox{10pt}[0pt][0pt]{\makebox[0pt]{\hspace{-28pt}$\scriptstyle\circ$}}}
\newcommand{\twocases}{Suppose
   \begin{list}{}{\leftmargin=40mm \labelwidth=20mm \labelsep=5mm}
    \item[{\hfill\rm \bf either}\hfill] $H$ is  strongly elliptic 
    \item[{\hfill\rm \bf or}\hfill] $G$ is stratified and 
          $a_1,\ldots,a_{d'}$ is a basis for $\gotg_1$ 
     \newline
   in the stratification $(\gotg_m)_{m \in \{ 1,\ldots,r \} }$ of $\gotg$.
        \end{list}    }

\newcommand{\lstar}[1]{\mathbin{{}_{#1}\mskip-1mu*}}
\newcommand{\ldoublestar}[2]{\mathbin{{}_{#2}\mskip-1mu({}_{#1}\mskip-1mu*)}}
\newcommand{\sdp}[1]{\mathbin{{}_{#1}\mskip-2mu\times}}
\newcommand{\tfrac}[2]{{\textstyle \frac{#1}{#2} }}
\newcommand{\sptilde}{\hspace{0pt} \widetilde{\rule{0pt}{8pt} \hspace{7pt} }}

\newcommand{\cc}{{\cal C}}
\newcommand{\cd}{{\cal D}}
\newcommand{\cf}{{\cal F}}
\newcommand{\ch}{{\cal H}}
\newcommand{\ck}{{\cal K}}
\newcommand{\cl}{{\cal L}}
\newcommand{\cm}{{\cal M}}
\newcommand{\co}{{\cal O}}
\newcommand{\cs}{{\cal S}}
\newcommand{\ct}{{\cal T}}
\newcommand{\cx}{{\cal X}}
\newcommand{\cy}{{\cal Y}}
\newcommand{\cz}{{\cal Z}}

\newfont{\fontcmrten}{cmr10}
\newcommand{\slbrl}{\mbox{\fontcmrten (}}
\newcommand{\slbrr}{\mbox{\fontcmrten )}}

\newlength{\hightcharacter}
\newlength{\widthcharacter}
\newcommand{\covsup}[1]{\settowidth{\widthcharacter}{$#1$}\addtolength{\widthcharacter}{-0.15em}\settoheight{\hightcharacter}{$#1$}\addtolength{\hightcharacter}{0.1ex}#1\raisebox{\hightcharacter}[0pt][0pt]{\makebox[0pt]{\hspace{-\widthcharacter}$\scriptstyle\circ$}}}
\newcommand{\cov}[1]{\settowidth{\widthcharacter}{$#1$}\addtolength{\widthcharacter}{-0.15em}\settoheight{\hightcharacter}{$#1$}\addtolength{\hightcharacter}{0.1ex}#1\raisebox{\hightcharacter}{\makebox[0pt]{\hspace{-\widthcharacter}$\scriptstyle\circ$}}}
\newcommand{\scov}[1]{\settowidth{\widthcharacter}{$#1$}\addtolength{\widthcharacter}{-0.15em}\settoheight{\hightcharacter}{$#1$}\addtolength{\hightcharacter}{0.1ex}#1\raisebox{0.7\hightcharacter}{\makebox[0pt]{\hspace{-\widthcharacter}$\scriptstyle\circ$}}}
\newcommand{\gznp}{\settowidth{\widthcharacter}{$G$}\addtolength{\widthcharacter}{-0.15em}\settoheight{\hightcharacter}{$G$}\addtolength{\hightcharacter}{0.1ex}G\raisebox{\hightcharacter}[0pt][0pt]{\makebox[0pt]{\hspace{-\widthcharacter}$\scriptstyle\circ$}}_N'}
\newcommand{\sgznp}{{\settowidth{\widthcharacter}{$G$}\addtolength{\widthcharacter}{-0.15em}\settoheight{\hightcharacter}{$G$}\addtolength{\hightcharacter}{0.1ex}\raisebox{-5pt}{$\scriptstyle G\raisebox{0.7\hightcharacter}[0pt][0pt]{\makebox[0pt]{\hspace{-\widthcharacter}$\scriptstyle\circ$}}_N'$}}}

 \thispagestyle{empty}

\begin{center}
{\Large{\bf Analysis of degenerate elliptic operators  }}\\[2mm] 
{\Large{\bf of  Gru\v{s}in type}}  \\[2mm]
\large  Derek W. Robinson$^1$ and Adam Sikora$^2$\\[2mm]

\normalsize{June 2006}
\end{center}

\vspace{5mm}

\begin{center}
{\bf Abstract}
\end{center}

\begin{list}{}{\leftmargin=1.8cm \rightmargin=1.8cm \listparindent=10mm 
   \parsep=0pt}
   \item
We analyze degenerate, second-order, elliptic operators $H$ in divergence form on $L_2(\Ri^{n}\times\Ri^{m})$.
We assume the coefficients are real symmetric and $a_1H_\delta\geq H\geq a_2H_\delta$ for some
$a_1,a_2>0$ where
\[
H_\delta=-\nabla_{x_1}\,c_{\delta_1, \delta'_1}(x_1)\,\nabla_{x_1}-c_{\delta_2, \delta'_2}(x_1)\,\nabla_{x_2}^2
\;\;\;.
\]
Here $x_1\in\Ri^n$, $x_2\in\Ri^m$ and $c_{\delta_i, \delta'_i}$ are positive measurable functions
such that $c_{\delta_i, \delta'_i}(x)$ behaves like $|x|^{\delta_i}$ as $x\to0$ and 
$|x|^{\delta_i'}$ as $x\to\infty$ with $\delta_1,\delta_1'\in[0,1\rangle$
 and  $\delta_2,\delta_2'\geq0$.

Our principal results state that the submarkovian semigroup $S_t=e^{-tH}$ is conservative and its
kernel $K_t$ satisfies  bounds
\[
0\leq  K_t(x\,;y)\leq a\,(|B(x\,;t^{1/2})|\,|B(y\,;t^{1/2})|)^{-1/2}
\]
where  $|B(x\,;r)|$ denotes the volume of the ball $B(x\,;r)$  centred at $x$ with radius $r$
measured with respect to the  Riemannian distance associated with $H$.
The proofs depend on detailed subelliptic estimations on $H$, a precise characterization of the 
Riemannian distance and the corresponding volumes and wave equation techniques which
exploit the finite speed of propagation.

We discuss further implications of these bounds and give explicit examples that show the kernel is not necessarily strictly positive, nor continuous.

\end{list}

\vspace{1cm}


\vspace{1cm}

\vspace{1cm}

\noindent
\begin{tabular}{@{}cl@{\hspace{10mm}}cl}
1. &Mathematical Sciences Institute  & 
  2. &Department of Mathematical Sciences  \\
& Australian National University  & 
  &   New Mexico State University  \\
& Canberra, ACT 0200& 
  &  P.O. Box 30001 \\
& Australia.& 
  & Las Cruces  \\
&{}& 
  & NM 88003-8001, USA.
\end{tabular}

\newpage
\setcounter{page}{1}

\section{Introduction}\label{Scsg1}

The classical work of Nash \cite{Nash}, De Giorgi \cite{DG} Aronson  \cite{Aro} gives boundedness and regularity properties of weak solutions of linear   parabolic and elliptic equations  with measurable coefficients.
Their  results cover the equations given by  second-order  operators $H$ in divergence form
on $L_2(\Ri^d)$, i.e., operators formally expressed as 
\begin{equation}
H=-\sum^d_{i,j=1}\partial_i\,c_{ij}\,\partial_j
\;\;\;,
\label{ecsg1.1}
\end{equation}
where $\partial_i=\partial/\partial x_i$, the $c_{ij}$ are real-valued measurable functions
and  the  coefficient matrix $C=(c_{ij})$ is assumed to be symmetric and  positive-definite 
almost-everywhere.
The principal  hypothesis of the theory is  (uniform)
strong ellipticity of the coefficients.
This condition is formulated as $\infty>\lambda I\geq C\geq \mu \,I>0$
or as  the equivalent operator condition $\lambda L\geq H\geq\mu\,L$ where $L=-\sum^d_{i=1}\partial_i^2$ is the usual Laplacian.
It is a condition of non-degeneracy.

The prime conclusion of the Nash--De Giorgi--Aronson theory is  that the fundamental solution of the parabolic equation, the heat kernel, satisfies global Gaussian upper and lower  bounds. 
H\"older continuity of solutions of the elliptic and parabolic equations  is  then a consequence of the Gaussian bounds.
It has become increasingly apparent that these  bounds encode a great deal of useful information on related concepts such as  Riesz transforms, spectral  multipliers, holomorphic functional calculus etc. (see, for example, \cite{DOS} \cite{Ouh5} and references therein).
It has also been established that the Gaussian bounds follow from two general structural features,
a Poincar\'e inequality and volume doubling \cite{Gri4} \cite{Sal4}.
This has led to extension of many 
results to a broader class of operators, e.g., the Laplace--Beltrami operator on Riemannian
manifolds with non-negative curvature (see \cite{Gri3} and \cite{Sal5} for  reviews of these developments).

The theory of  degenerate elliptic operators is in comparison underdeveloped although two general classes have been identified which retain many of the structural
features associated with strong ellipticity, at least locally, albeit with significant geometric modifications.
These classes weaken the strong ellipticity hypotheses in different directions and 
cover complementary types of degeneracy but only describe some of the possibilities .

The first class   consists of operators for which the largest eigenvalue $\mu_M$ and the  inverse of the smallest eigenvalue $\mu_m$ of the coefficient
matrix $C$ are both  locally integrable  \cite{FKS} \cite{Tru2}.
 Poincar\'e and Harnack inequalities, H\"older continuity etc. follow from  this condition, or slightly more stringent conditions.
 These conditions place direct restraints on the order of degeneracy, e.g.,  the order of $\mu_m$
 in the neighbourhood of a zero, and limit the analysis to weakly degenerate operators.
 Under such restraints one can still exploit many of the methods used to analyze strongly elliptic   operators.
Note, however, that this class of operators does not contain simple examples such as the Heisenberg sublaplacian $H_{\rm Heis}=-\partial_1^2-(\partial_2+x_1\,\partial_3)^2$ on $L_2(\Ri^3)$ for which $\mu_m$ 
is identically zero.

The second class consists of operators which satisfy a subelliptic estimate
$H\geq\mu\,L^\gamma-\nu\,I$ where $\mu>0$, $\nu\geq0$ and $\gamma\in\langle0,1]$.
This condition, which  is a natural extension of  the operator form of strong ellipticity $H\geq\mu\,L$,
 is of a rather different nature  since it cannot be expressed in terms of 
the  lowest eigenvalue $\mu_m$  or, indeed,  in terms of  the coefficients~$c_{ij}$.
For example, the Heisenberg sublaplacian $H_{\rm Heis}$ on $L_2(\Ri^3)$ satisfies the estimate, with $\gamma=1/2$. 
This is a specific example of a type of subelliptic operator introduced by H\"ormander \cite{Hor1}.
The H\"ormander operators are expressed in the form $H=\sum^n_{i=1}X_i^*X_i$ where the $X_i$
are smooth vector fields.
The principal assumption is that the $X_i$ together with their multicommutators of a fixed order
$r$ span the tangent space at each point $x\in\Ri^d$.
Then $H$ satisfies the subellipticity condition locally with $\gamma=1/r$ (see, \cite{JSC1} for
a detailed review  and references).
The  algebraic structure of the fields required by the H\"ormander commutator condition
places  restrictions on the smoothness and growth of  the coefficients of the operators.
Fefferman and Phong \cite{FP} subsequently extended the scope of the theory by establishing
that  locally the subellipticity is equivalent to a   property of the balls 
defined by the corresponding Riemannian geometry.
The results of these authors then opened the way to the extension of many of the local estimates
of the Nash--De Giorgi--Aronson theory to the subelliptic situation (again see, \cite{JSC1} for
a detailed description). 
Nevertheless the results of Fefferman and Phong depend on some smoothness of the coefficients
of $H$.
In particular it is necessary for the coefficients to be at least twice differentiable to define the
balls that play a key role in their analysis (see \cite{FSC}, Section~1) and this is possibly sufficient.
This is indeed the case  in two dimensions Xu \cite{Xu1}.
But the  $C^2$ requirement  places a substantial  restraint on the possible degeneracy.

In this paper we introduce and analyze global properties of  a family of degenerate operators with measurable coefficients which is not covered  by either of the above classes but which incorporates  two different types of degeneracy typical of the  classes.
As a preliminary we  define two functions $f$, $g$ with values in an ordered space  to be equivalent, $f\sim g$, 
if there are $a$, $a'>0$ such that $a\,f\leq g\leq a'\,f$ uniformly.
We use this notion in a variety of contexts, e.g.,  for functions over $\Ri^d$, for positive symmetric matrices and  lower semibounded selfadjoint operators and 
also for  the equivalence of  quadratic forms.
In addition  for $a>0$ we set
\begin{eqnarray*}
  a^{(\alpha,\alpha')}=   \left\{ \begin{array}{llll}
 a^\alpha
 & \mbox{ if $ a \leq 1 $}\\[5pt]
  a^{\alpha'}
&     \mbox{ if
$   a \geq 1  $}
           \end{array}
    \right.
\end{eqnarray*}
for all $\alpha, \alpha'\in\Ri$.

We analyze   operators $H$  on $L_2(\Ri^n\times \Ri^m)$  of the form (\ref{ecsg1.1})
but with the matrix of coefficients  $C\sim C_\delta$
where  $C_\delta$ is  a block diagonal matrix
$C_\delta(x_1,x_2)=c_{\delta_1, \delta'_1}(x_1)\, I_n+c_{\delta_2, \delta'_2}(x_1)\,I_m$
and $c_{\delta_1, \delta_1'},c_{\delta_2, \delta'_2}$ are  positive measurable  functions such that
\begin{equation}
c_{\delta_i, \delta_i'}(x)\sim |x|^{(2\delta_i,2\delta_i')}
\label{ecsg1.10}
\end{equation}
where
the indices $\delta_1,\delta_2,\delta'_1,\delta'_2$ are all non-negative and $\delta_1,\delta_1'<1$
but there is no upper bound on $\delta_2$ and $\delta'_2$.
Thus  $H\sim H_\delta$,
where $H_\delta$ is given by 
\begin{equation}
H_\delta=-\nabla_{x_1}\,c_{\delta_1, \delta'_1}\,\nabla_{x_1}+c_{\delta_2, \delta'_2}\,L_{x_2}
\label{ecsg1.2}
\end{equation}
with  $x_1\in\Ri^n$, $x_2\in\Ri^m$, $\nabla_{x_1}$  the gradient operator on $L_2(\Ri^n)$ and  $L_{x_2}=-\nabla_{x_2}^2$  the Laplacian on $L_2(\Ri^m)$.
Since $C_\delta$ is only defined up to equivalence there is a freedom of choice which will often be
exploited in the sequel.

In summary we consider elliptic operators $H$ of the form (\ref{ecsg1.1}), defined precisely by quadratic
form techniques in Section~\ref{Scsg2}, with $H\sim H_\delta$ where $H_\delta$ is the elliptic operator 
(\ref{ecsg1.2}) with coefficients  satisfying (\ref{ecsg1.10}).
The operator (\ref{ecsg1.2}) is a natural generalization of the H\"ormander type  operators  
$-\partial_1^2-x_1^{2n}\,\partial_2^2$ introduced by Gru\v{s}in \cite{Gru}.
Therefore we refer to them as Gru\v{s}in operators.
Note that for these operators the lowest eigenvalue $\mu_m$ of the coefficient matrix $C$ satisfies 
$\mu_m(x)\sim |x_1|^{2(\delta_1\vee\delta_2)}$ for $|x|\leq 1$ thus the inverse of $\mu_m$ is locally integrable if 
and only if $\delta_1\vee\delta_2<n$.
Moreover,  we make no smoothness assumptions on the coefficients $c_{ij}$.

We will prove that the semigroup $S$ generated by the  Gru\v{s}in operator $H$ is conservative, i.e., $S_t\one=\one$ on $L_\infty(\Ri^{n+m})$ and that its kernel $K$ satisfies Gaussian bounds with respect to the appropriate Riemannian geometry.

The Riemannian distance associated with an elliptic operator 
 (\ref{ecsg1.1}) with measurable coefficients $C=(c_{ij})$ is defined by 
 \begin{equation}
d(y\,;z)=\sup_{\psi\in D}(\psi(y)-\psi(z))
\label{ecsg1.21}
\end{equation}
for all  $y,z\in \Ri^{d}$ where 
\begin{equation}
D=\{\psi\in W^{1,\infty}(\Ri^d):\,
\sum^d_{i,j=1}c_{ij}(\partial_i\psi)(\partial_j\psi)\leq 1\,\}
\;\;\;.
\label{ecsg1.211}
\end{equation}
If $C$ is strongly elliptic then this is equivalent to the usual Euclidean distance but for a general
elliptic operator it is not necessarily a genuine distance, i.e., it could take the value  infinity. 
If, however, $H$ is a  Gru\v{s}in operator on $\Ri^n\times \Ri^m$ then $C\sim C_\delta$ and 
$d(\cdot\,;\cdot)$ is equivalent to 
\[
d_\delta(y\,;z)=\sup\{|\psi(y)-\psi(z)|:\psi\in W^{1,\infty}(\Ri^n\times \Ri^m)\,,\,
c_{\delta_1, \delta_1'}\,|\nabla_{x_1}\psi|^2+c_{\delta_2, \delta'_2}\,|\nabla_{x_2}\psi|^2\leq 1\,\}
\;\;\;.
\]
The restrictions  $\delta_1,\delta_1'<1$ ensure that this latter function is a {\it bona fide} distance.

There are two technical difficulties in our analysis of Gru\v{s}in operators.
First one needs to derive subelliptic estimates on $H_\delta$ in order to obtain {\it a priori}
estimates on the semigroup $S$  (see Section~\ref{Scsg3}).
Secondly it is necessary to estimate
the volume (Lebesgue measure)  $|B(x\,;r)|$ of the Riemannian balls $B(x\,;r)=\{y:d(x\,;y)<r\}$ or,
equivalently, the balls  $B_\delta(x\,;r)=\{y:d_\delta(x\,;y)<r\}$ (see Section~\ref{Scsg5}).
But then we establish in Section~\ref{Scsg6} that the semigroup $S$ is conservative and and its kernel
$K$ satisfies Gaussian bounds.
Specifically, for each $\varepsilon\in\langle0,1]$ there is an $a>0$ such that
\[
0\leq K_t(x\,;y)\leq a\,(|B(x\,;t^{1/2})|\,|B(y\,;t^{1/2})|)^{-1/2}\,e^{-d(x;y)^2/(4(1+\varepsilon)t)}
\]
for all $x,y\in\Ri^n\times\Ri^m$ and all $t>0$.

Lower bounds and continuity properties are more sensitive, e.g., if $n=1$ and $\delta_1\in[1/2,1\rangle$
then the kernel is not strictly positive nor is it continuous  (see \cite{ERSZ1}, Sections~5 and 6)
since the system separates into two distinct subsystems, the halfspace $x_1\geq0$ and the halfspace
$x_1\leq0$.
In Section~\ref{Scsg7} we discuss  the simplest example in one-dimension, $n=1$, $m=0$.
The one-dimensional example has been analyzed at length in the setting of control theory
\cite{ABCF} \cite{MaV} from a slightly different perspective.
In the latter context the case $\delta\in[0,1/2\rangle$ is referred to as weakly degenerate
and the case $\delta_1\in[1/2,1\rangle$ as strongly degenerate.
We will indeed establish positivity and continuity properties similar
to those associated with strongly elliptic operators  in  the weakly degenerate case.
In addition we establish some partial regularity and positivity properties 
for the ergodic components on $L_2(\Ri_\pm)$ in the strongly degenerate case.
We emphasize that the separation phenomenon
demonstrates that the Riemannian distance is not appropriate for the full
description of the structure associated with strongly  degenerate operators.
Indeed the small time asymptotics of the kernel is given by a larger distance \cite{ERS4} \cite{HiR}
\cite{AH} which incorporates the separation.
This indicates that the above Gaussian bounds are not optimal.
Nevertheless, they suffice for the derivation of several significant results in the global spectral analysis
of  Gru\v{s}in operators (see Section~\ref{Scsg8}).

Finally we note that Sawyer and Wheeden \cite{SWh} have recently reformulated the 
H\"ormander--Fefferman--Phong theory to incorporate elliptic operators with $L_\infty$-coefficients.
But their emphasis is from the outset on establishing local regularity properties of weak solutions.
Our analysis of Gru\v{s}in operators with  $n=1$ and $\delta_1\in[1/2,1\rangle$
shows that such properties are not universally valid.


\section{Preliminaries}\label{Scsg2}

In the sequel we are only interested in qualitative estimates.
Hence  we adopt the variable constant convention.
In subsequent estimates  $a$, $a'$ etc.  denote strictly positive constants whose value may vary
bound by bound
but which are independent of the key variables in the estimates.

We begin by discussing  the precise definition of  elliptic operators of the form (\ref{ecsg1.1}) or more specifically (\ref{ecsg1.2}) through closed quadratic forms.

First, let $C=(c_{ij})$ be a real-valued symmetric matrix with measurable, locally integrable, coefficients $c_{ij}$.
Assume $C$ is positive-definite almost everywhere and  
define  the   positive quadratic  form $h$  by
\begin{equation}
h(\varphi)=\sum^d_{i,j=1}(\partial_i\varphi,c_{ij}\partial_j\varphi)
= \int_{\Ri^d}dx\,\sum^d_{i,j=1}c_{ij}(x)(\partial_i\varphi)(x)(\partial_j\varphi)(x)
\label{ecsg2.1}
\end{equation}
with $D(h)$ consisting of those $\varphi\in W^{1,2}(\Ri^d)$ for which the integral converges.
Since the coefficients are locally integrable the form is densely defined.
But it is not necessarily closed.
There are three possibilities.

First, if $h$ is closed then there is a positive self-adjoint operator $H$  such that $D(h)=D(H^{1/2})$ and
 $h(\varphi)=\|H^{1/2}\varphi\|_2^2$.
Then  we define $H$ to be the elliptic operator with coefficients $C$.
For example if $C$ satisfies the usual strong ellipticity assumptions, $\lambda I\geq C \geq \mu I>0$, then $D(h)=W^{1,2}(\Ri^d)$ and $h$ is closed.
More generally one has the following.

\begin{lemma}\label{lcsg1.1}
 If $C\geq \mu I>0$ almost everywhere then $h$ is closed. 
 \end{lemma}
\proof\
The lemma is established by a  monotonicity argument.
If  $C_N=C\wedge NI$ with $N\in \Ni$ then $C_N$ is strongly elliptic and  the   corresponding closed forms $h_N$ with $ D(h_N)=W^{1,2}(\Ri^d)$ are monotonically increasing.
But  $h(\varphi)=\sup_Nh_N(\varphi)$ with  $D(h)$ the subspace of 
$W^{1,2}(\Ri^d)$ for which the supremum is finite.
Then, however, $h$ is closed (for details see, for example \cite{BR2}, Lemma 5.2.23).
 \hfill$\Box$

\begin{remarkn}\label{rcsg1.1} Note that if $H$ is the elliptic operator associated with the closed form $h$ in the lemma and $H_N$ the strongly elliptic operators associated with the $h_N$ then the $H_N$ converge in the strong resolvent sense to $H$.
But the $h_N$ are Dirichlet forms.
Hence $h$ is a Dirichlet form and the semigroup $S$ generated by $H$ is submarkovian.
\end{remarkn}

Secondly, it is possible that $h$ is closable although it is not closed.
Then the elliptic operator $H$ is naturally defined as the self-adjoint operator associated with the 
closure $\overline h$ of $h$. 
Again $\overline h$  is a Dirichlet form and the elliptic semigroup $S$ is submarkovian.
This is indeed the case of principal interest in the current context.

\begin{lemma}\label{lcsg1.2}
The form $h$ of a Gru\v{s}in operator is closable.
\end{lemma}
\proof\
Let $h_\delta$ denote the form of the operator $H_\delta$ given formally by (\ref{ecsg1.2}).
Since $C\sim C_\delta$ it follows that $h\sim h_\delta$.
Explicitly, $D(h)=D(h_\delta)$ and there are $a,a'>0$ such that 
\[
a\,h_\delta(\varphi)\geq h(\varphi)\geq a'\,h_\delta(\varphi)
\]
for all $\varphi\in D(h)$.
Thus it suffices to prove that $h_\delta$ is closable.
But we may assume the coefficients  $c_{\delta_1, \delta_1'},c_{\delta_2, \delta'_2}$ are continuous.
Then $h_\delta$ is closable by the proof of Proposition~2.3 in \cite{ERSZ1}.
\hfill$\Box$

\bigskip

Although the  forms of the Gru\v{s}in operators are closable it is nevertheless of interest to 
consider the third possibility for general elliptic forms, the possibility that 
 $h$  is neither closed nor  closable.
(For explicit examples see  \cite{FOT} Section~3.1 and in particular Theorem~3.1).
In this case   one can introduce the relaxation $h_0$ of $h$ as the  largest  positive, closed, quadratic form $h_0$  such that $h_0\leq h$.
The relaxation occurs in the context of nonlinear phenomena and
discontinuous media  (see, for example,  \cite{Bra} \cite {ET} \cite{Jos} \cite{DalM}  \cite{Mosco}
and references therein) and can be characterized in several different ways.
 Simon~\cite{bSim5}, Theorems~2.1 and 2.2, has shown that  a general positive 
quadratic form $h$ 
 can be decomposed as a sum $h=h_r+h_s$ of two positive forms with 
$D(h_r)=D(h)=D(h_s)$ with
$h_r$  the largest closable form majorized by $h$.
Simon refers to $h_r$ as the regular part of $h$.
Then $h_0=\overline h_r$.
Note that  $h_0=h$  if $h$ is closed and  $h_0=\overline h$ if $h$ is closable.

There is an alternative method of constructing $h_0$ by monotone approximation.
Let $l$ be the closed quadratic form associated with the Laplacian $L$.
Then $D(l)=W^{1,2}(\Ri^d)$.
Next define  $h_\varepsilon=h+\varepsilon \,l$ for $\varepsilon>0$ with 
$D(h_\varepsilon)=D(h)$.
Then $h_\varepsilon$ corresponds to the elliptic operator with coefficients
$c_{ij}+\varepsilon\, \delta_{ij}$ and  is closed by Lemma~\ref{lcsg1.1}.
The corresponding positive self-adjoint operators $H_\varepsilon$ form a decreasing sequence which, by a result of Kato \cite{Kat1}, Theorem~VIII.3.11, converges in the strong
resolvent sense to a positive self-adjoint operator $H_0$.
This is 
 the operator $H_0$ associated with the relaxation $h_0$  (see \cite{bSim5}, Theorem~3.2).
The latter construction justifies the interpretation of $H_0$ as the elliptic operator with coefficients $C$.
It follows again that the relaxation $h_0$ is a Dirichlet form. 
Moreover the construction of the relaxation respects order properties, i.e., if $h$ and $k$ are two 
elliptic forms and $h\geq k$ then $h_0\geq k_0$.

In each of the above situations the elliptic operator $H$  with coefficients $C$  is obtained as a double limit $N\to\infty$, $\varepsilon\to0$ of the strongly elliptic operators $H_{N,\varepsilon}$ with coefficients $C_{N,\varepsilon}=(C\wedge NI)+\varepsilon I$.
The convergence is in the strong resolvent sense.
In particular this construction is applicable to the Gru\v{s}in operators.
Therefore it is not essential to make a notational distinction between the three cases.
Thus in the sequel we use $H$ to denote the self-adjoint elliptic operator constructed by this limiting process,  $S$ the submarkovian semigroup generated by $H$ and $K$ the semigroup kernel.
It follows in particular from the construction  that the semigroups $S^{(N,\varepsilon)}$ generated by the strongly elliptic operators  $H_{N,\varepsilon}$ 
converge strongly to $S$ on $L_2(\Ri^d)$.

 Our next  aim is to discuss {\it a priori} bounds on $t\to K_t$ which are uniform over $\Ri^d$.
 This is equivalent to obtaining bounds  on the crossnorms $\|S_t\|_{1\to\infty}$ of the semigroup 
 as a map from $L_1$ to $L_\infty$, or the crossnorms $\|S_t\|_{2\to\infty}$ from $L_2$ to $L_\infty$
 since 
   \[
 \|S_t\|_{1\to\infty}=\esssup_{x, y\in\Ri^d}|K_t(x\,;y)|
 =\esssup_{x\in\Ri^d}\int_{\Ri^d}dy\,|K_{t/2}(x\,;y)|^2=(\|S_{t/2}\|_{2\to\infty})^2
 \]
 for all $t>0$.

The standard method of obtaining bounds on the crossnorms $\|S_t\|_{2\to\infty}$, and hence on 
$\|S_t\|_{1\to\infty}$, for strongly elliptic semigroups is via Nash inequalities.
In the context of degenerate elliptic operators it is useful to consider a particular class of inequalities 
defined in terms of Fourier  multipliers.

Let $ F$ be a positive real function over $\Ri^d$ and define the corresponding Fourier multiplier  $F$  on 
$L_2(\Ri^d)$ by 
$\widetilde{(F\varphi)}(p)= F(p)\widetilde\varphi(p)$
where  $\widetilde\varphi$ denotes
the Fourier transform of $\varphi\in L_2(\Ri^d)$.
One can also interpret $F$ as  a differential operator $F=F(i\nabla_x)$ with constant coefficients.
Next  let $f $ denote the closed form  corresponding to $F $, i.e., 
\[
f(\varphi)=\int_{\Ri^d}dp\,F(p)\,|\widetilde\varphi(p)|^2
\]
with $D(f)$ the subspace of $\varphi\in L_2(\Ri^d)$ for which the integral is finite.
Finally  let $V_F(r)$ denote the volume (Lebesgue measure) of the set $\{p: F(p)<r^2\}$.

A subelliptic estimate of the form $h\geq f$ immediately gives a Nash type inequality.

\begin{lemma} \label{lcsg2.1}
If  $h\geq f$ then
\begin{equation}
\|\varphi\|_2^2\leq r^{-2}h(\varphi)
+(2\pi)^{-d}\,V_F(r)\,\|\varphi\|_1^2
\label{epre1.7}
\end{equation}
for all $\varphi\in D(h)\cap L_1$ and all $r>0$.
\end{lemma}
\proof\
The proof is a direct consequence of 
 the Plancherel formula;
\begin{eqnarray}
\|\varphi\|_2^2=\int_{\Ri^d}dp\,|\widetilde\varphi(p)|^2&=&
\int_{ F(p)\geq r^2}dp\,|\widetilde\varphi(p)|^2
+
\int_{ F(p)<r^2}dp\,|\widetilde\varphi(p)|^2\nonumber\\[5pt]
&\leq &
r^{-2}\int_{\Ri^d}dp\, F(p)|\widetilde\varphi(p)|^2
+
(2\pi)^{-d}\,\int_{ F( p)<r^2}dp\,\|\varphi\|_1^2\nonumber\\[5pt]
&\leq & r^{-2}h(\varphi)
+ (2\pi)^{-d}\,V_F(r)\,\|\varphi\|_1^2\label{epre1.71}
\end{eqnarray}
for all $r>0$.
\hfill$\Box$

\bigskip

The  Nash inequality allows one to obtain bounds on the cross-norm $\|S_t\|_{2\to\infty}$ for many different $F$. 
In particular if $V_F$ has a polynomial behaviour one can estimate $\|S_t\|_{1\to 2}$,
and by duality  $\|S_t\|_{2\to\infty}$, by a straightforward extension of Nash's original argument.
In particular one obtains  the following.

\begin{lemma}\label{lcsg2.2}
If $h\geq f$  and $V_F(r)\leq a\,r^{(D', D)}$ then 
$\|S_t\|_{1\to\infty}\leq b\,t^{(-D/2,-D'/2)}$.
\end{lemma}
\proof\ The result can be deduced from \cite{CKS}, Theorem~2.9 or from the alternative argument given in \cite{Robm} pages 268--269.\hfill$\Box$

\bigskip

The lemma demonstrates that  large values of $r$ give small $t$ bounds and small values of~$r$ give large $t$ bounds.
If $h$ is strongly elliptic then $ F(p)=\mu\, p^2=(\mu^{1/2}p)^2$, $V_F(r)\sim r^d$ and  $\|S_t\|_{1\to\infty}\leq a\,(\mu t)^{-d/2}$ for all $t>0$.
Alternatively the subellipticity condition $H\geq\mu\,L^\gamma-\nu\,I$ corresponds to $ F(p)=(\mu\,|p|^{2\gamma}-\nu)\vee 0$ and this  only gives useful information on the large~$r$ behaviour of $V_F$.
It  yields bounds
$\|S_t\|_{1\to\infty}\leq a\,t^{-d/(2\gamma)}$ for $t\leq 1$.

\smallskip

Lemmas~\ref{lcsg2.1} and \ref{lcsg2.2}  allow us in Section~\ref{Scsg3} to obtain uniform bounds on the semigroup kernels associated with the  Gru\v{s}in  operators (\ref{ecsg1.2}) except if $n=1$ and $\delta_1$ or $\delta_2$  is in
$[1/2,1\rangle$.
In the latter case one obtains subelliptic bounds of a different character, bounds in terms of the Neumann Laplacian.
But these can also be used to obtain Nash inequalities.

   Let $ L_{x,N}$ denote the self-adjoint version of  the operator $-d^2/dx^2$ on $L_2(\Ri)$ with Neumann
   boundary conditions at the origin.
   Then $L_{x,N}=L^+_{x,N}\oplus L^+_{x,N}$ where $L^{\pm}_{x,N}$ are the Neumann operators on
   $L_2(\Ri_\pm)$, respectively.
   Now let $F$ be a positive function on the half-line and define the operators  $F_N=F(L_{x,N})$ and
  $F^\pm_N=F(L^\pm_{x,N})$   by spectral theory. 
   Then $F_N=F_N^+\oplus F^-_N$.
   Next let $f_N$ and $f^\pm_N$ be the  forms corresponding to $F_N$ and $F_N^\pm$.
   Each $\varphi\in D(f_N)$ has a unique decomposition $\varphi=\varphi_+\oplus\varphi_-$
   with $\varphi_\pm\in D(f_N^\pm)$ and one has the relations $\|\varphi\|_2^2=\|\varphi_+\|_2^2+\|\varphi_-\|_2^2$,
   $f_N(\varphi)=f^+_N(\varphi_+)+f^-_N(\varphi_-)$ and $\|\varphi\|_1=\|\varphi_+\|_1+\|\varphi_-\|_1$.
   Thus the Nash inequalities can be analyzed by examining the two subsystems on $L_2(\Ri_\pm)$.
   This can then be handled by considering the extension of the operators from 
$L_2(\Ri_\pm)$ to the space $L_{2,e}(\Ri)$ of even  functions over the line.
Since both cases are similar we only consider the extension of $L_2(\Ri_+)$.

First, if $\varphi\in L_2(\Ri_+)$ we define the symmetric extension $E\varphi$ to 
$L_2(\Ri)$ by $(E\varphi)(\pm x)=\varphi(x)$ if $x \geq0$.
Then $L_{2,e}(\Ri)$  is the range of $E$.
Secondly, if $\varphi\in L_2(\Ri)$ we define the restriction $R\varphi$ to $L_2(\Ri_+)$ by
$(R\varphi)(x)=\varphi(x)$ if $x\geq 0$.
Thirdly, if $T$ denotes the semigroup generated by the one-dimensional
Laplacian $L_{x}=-d^2/dx^2$  on $L_2(\Ri)$ then $T$
leaves  $L_{2,e}(\Ri)$  invariant.
Hence $RTE$ defines a semigroup on $L_2(\Ri_+)$. 
The semigroup property follows because 
\[
(RT_sE)(RT_tE)=(RT_s)(ER)(T_tE)=(RT_s)(T_tE)=(RT_{s+t}E)
\]
where the second step uses the invariance.
Now it is easy to check that $RT_tE=T^{(N)}_t$ where $T^{(N)}$ denotes 
the semigroup generated by $L_{x,N}$.
Similar considerations apply to $L_{x_1,N}\otimes \one$ on $L_2(\Ri\times\Ri^m)$.
If $E\varphi$ is the even  extension of $\varphi\in L_2(\Ri_+\times\Ri^m)$ to $L_2(\Ri\times\Ri^m)$
and $R$ the corresponding restriction then one has $R(T_t\otimes\one)E=T^{(N)}_t\otimes\one$. 

  Next let $F$ be a positive bounded function on $\Ri_+\times \Ri^m$.
  Define the multiplier by $F$ by 
  $\widetilde{(F\varphi)}(p_1,p_2)= F(p_1^2,p_2)\widetilde\varphi(p_1,p_2)$
and the corresponding operators on $L_2(\Ri\times\Ri^m)$ 
 by  $F=F(L_{x_1},i\nabla_{x_2})$  and $F_N= F(L_{x_1,N},i\nabla_{x_2})$.
 Let $F_N^+$ denote the restriction of $F_N$ to  $L_2(\Ri_+\times\Ri^m)$.
Further let  $f$, $f_N$  and $f^+_N$ be the corresponding forms. 
Then if $\varphi\in D(f_N^+)$ one has $E\varphi\in D(f)$ and 
\begin{eqnarray*}
f(E\varphi)=(E\varphi, FE\varphi)&=&(E\varphi, ERFE\varphi)
=(E\varphi, EF_N^+\varphi)=2\,(\varphi, F_N^+\varphi)=2\,f_N^+(\varphi)
\;\;\;.
\end{eqnarray*}
Hence using the estimate (\ref{epre1.71}) for $F$ one deduces that 
\begin{eqnarray*}
\|\varphi\|_2^2=2^{-1}\|E\varphi\|_2^2
&\leq &2^{-1}r^{-2}(E\varphi, FE\varphi)+2^{-1}V_F(r) \,\|E\varphi\|_1^2
=r^{-2}f_N^+(\varphi)+2\,(2\pi)^{-d}\,V_F(r) \,\|\varphi\|_1^2
\end{eqnarray*}
for all $r>0$.

One can reason similarly for the restriction to the left half line and then by combination one obtains the 
following analogue of Lemma~\ref{lcsg2.1} with $V_F(r)$  the volume of the set $\{(p_1,p_2): F(p_1^2, p_2)<r^2\}$.

\begin{lemma} \label{lcsg2.3}
If  $h\geq f_N$ in the form sense on $L_2(\Ri\times\Ri^m)$  then
\begin{equation}
\|\varphi\|_2^2\leq r^{-2}h(\varphi)
+ 4\,(2\pi)^{-d}\,V_F(r)\,\|\varphi\|_1^2
\label{epre1.8}
\end{equation}
for all $\varphi\in D(h)\cap L_1(\Ri\times\Ri^m)$ and all $r>0$.
\end{lemma}

If the function $V_F$ in the Nash inequality of Lemma~\ref{lcsg2.3} has a polynomial growth of the type considered
in Lemma~\ref{lcsg2.2} one then obtains analogous bounds on the crossnorm of the semigroup $S_t$.
Thus the key point  in the analysis of the Gru\v{s}in operator with $n=1$ is to obtain a subelliptic estimate $h\geq f_N$ with the modified Fourier multiplier involving the Neumann Laplacian in the first direction.

\section{Subelliptic estimates} \label{Scsg3}

In this section we examine the Gru\v{s}in  operators and  derive uniform estimates
on the semigroup crossnorms $\|S_t\|_{1\to\infty}$ by use of the Nash inequalities of 
Lemmas~\ref{lcsg2.1} and \ref{lcsg2.3}.
These bounds will then be improved by other techniques in the sequel.
The principal result of this section is given by the following proposition.

\begin{prop}\label{pcsg3.0}
Let $S_t$ denote the positive self-adjoint semigroup on $L_2(\Ri^n\times\Ri^m)$ generated
by the Gru\v{s}in  operator  $H$  with coefficients $C\sim C_\delta$ where $C_\delta$  satisfies $(\ref{ecsg1.10})$.
Then
\[
\|S_t\|_{1\to\infty}\leq a\,t^{(-D/2,-D'/2)}
\]
 where
\[
D=(n+m(1+\delta_2-\delta_1))(1-\delta_1)^{-1}\;\;\;\;\;{and}\;\;\;\;\;
D'=(n+m(1+\delta_2'-\delta_1'))(1-\delta_1')^{-1}
\;\;\;.
\]
\end{prop}

Note that the local dimension $D$ depends only on the parameters $\delta_1$ and $\delta_2$ which govern
the local degeneracies of the coefficients of $H$ and the global dimension $D'$ depends on the parameters
$\delta_1'$ and $\delta_2'$ which govern
the global  degeneracies.
Moreover, $D,D'\geq n+m$, the Euclidean dimension, with $D=n+m$ if and only if $\delta_1=0=\delta_2$
and $D'=n+m$ if and only if $\delta_1'=0=\delta_2'$.

The proof of the proposition is in two stages.
First consider the operator
\begin{equation}
H_1=-\nabla_x\,c_{\delta_1, \delta_1'}\,\nabla_x
\label{esub1}
\end{equation}
 on $L_2(\Ri^n)$. 
There are two possibilities.

\medskip

\noindent {\bf Either} $\delta_1\geq \delta_1'$ and then $c_{\delta_1, \delta_1'}(x)\sim 
|x|^{2\delta_1}(1+|x|^{2\delta_1})^{-1+\delta_1'/\delta_1}$,

\medskip

\noindent {\bf or} $\hspace{9mm}\delta_1\leq \delta_1'$ and then $c_{\delta_1, \delta_1'}(x)
\sim |x|^{2\delta_1}+|x|^{2\delta_1'}$.

\medskip

\noindent Since the subsequent estimates are valid up to equivalence we can effectively replace $c_{\delta_1,\delta_1'}$
by the appropriate  function on the right.
Moreover, the form
\[
h_1(\varphi)=(\nabla_x\varphi,c_{\delta_1, \delta_1'}\,\nabla_x\varphi)
=\int_{\Ri^n}dx\,c_{\delta_1, \delta_1'}(x)|(\nabla_x\varphi)(x)|^2
\]
is closed on the domain $D(h_1)=W^{1,2}(\Ri^n\,;c_{\delta_1, \delta_1'}dx)$
and $C_c^\infty(\Ri^n)$ is a core of $h$.
Therefore it suffices to establish the following form estimates on $C_c^\infty(\Ri^n)$.

It is convenient to express the subelliptic estimates in operator terms.
Recall that $L_x=-\nabla_x^2$ denotes the Laplacian.

\begin{prop}\label{pcsg3.1}
Let $n\geq2$, or $n=1$ and $\delta_1, \delta_1'\in[0,1/2\rangle$.
Then $h_1\geq f$ on $L_2(\Ri^n\times\Ri^m)$ where $f$ is the form of the operator $F$  with
\[
 F\sim L_x^{(1-\delta_1')}(1+ L_x)^{-(\delta_1-\delta_1')}
\]
if $\delta_1\geq\delta_1'$ and 
\[
 F\sim  L_x^{(1-\delta_1)}+  L_x^{(1-\delta_1')}
\]
if $\delta_1\leq\delta_1'$.

Moreover, if   $n=1$ and $\delta_1\in[1/2,1\rangle$ or $\delta_1'\in[1/2,1\rangle$ then $h_1\geq f_N$  where $f_N$ is the form of the operator $F_N$  obtained by replacing $L_x$ by $L_{x,N}$.
\end{prop}

The two key properties needed for subelliptic estimates of $H_1$ are the following. 

\begin{lemma}\label{lsub1}
If $\gamma\in[0, 1\wedge n/2\rangle$ then
\[
L_x^{\,\gamma}\geq a\,|x|^{-2\gamma}
\]
in the form sense on $L_2(\Ri^n)$.

Moreover, if $n=1$ then $L_{x,D}\geq (4x^2)^{-1}$ on $L_2(\Ri)$ where $L_{x,D}$ is the Laplacian
with Dirichlet boundary conditions at the origin.
\end{lemma}

The statements are   versions of Hardy's inequality (see, for example, \cite{Dav13} and references
therein) and special cases of the inequalities of Caffarelli,  Kohn and Nirenberg \cite{CKN}.
The multidimensional version is often stated with $\gamma=1$ and $n\geq 3$.
In the latter case one has $a=(n-2)^2/4$ and this value is optimal. 
The fractional version follows from Strichartz' work \cite{Stri} on Fourier multipliers.

\begin{lemma}\label{lsub2}
Let $A$ and $B$ be self-adjoint operators on $L_2(\Ri^n)$ and let $\gamma\in[0, 1]$.
If  $A \geq B\geq0$ then
\[
A(I+A)^{-\gamma}\geq B(I+B)^{-\gamma}
\]
in the form sense.
\end{lemma}
\proof\
First one has
\[
A(\lambda I+A)^{-1}=I-\lambda\,(\lambda I+A)^{-1}\geq I-\lambda\,(\lambda I+B)^{-1}
=B(\lambda I+B)^{-1}
\]
for all $\lambda>0$ and this gives the result for $\gamma=1$.
But if $\gamma<1$ then
   \begin{eqnarray*}
  A(I+A)^{-\gamma}&=&c_\gamma\int^\infty_0{{d\lambda}\over{\lambda^{\gamma}}}\,A((1+\lambda)I+A)^{-1}
  \\[5pt]
  &\geq& c_\gamma\int^\infty_0{{d\lambda}\over{\lambda^{\gamma}}}\,B((1+\lambda)I+B)^{-1}
  =B(I+B)^{-\gamma}
  \end{eqnarray*}
  where we have used the standard integral representation of the fractional power.
  \hfill$\Box$
  
  \bigskip
 
 \noindent{\bf Proof of Proposition~\ref{pcsg3.1}}
Consider the case $\delta_1\geq \delta_1'$.
 Then for $\varphi\in C_c^\infty(\Ri^n)$ 
 \begin{eqnarray*}
 h_1(\varphi)&\geq&a\,(\nabla_x\varphi,|x|^{2\delta_1}(1+|x|^{2\delta_1})^{-1+\delta_1'/\delta_1}\nabla_x\varphi)\\[5pt]
 &\geq&a\,(\nabla_x\varphi,L_x^{-\delta_1}(1+L_x^{-\delta_1})^{-1+\delta_1'/\delta_1}\nabla_x\varphi)
 \geq a\,(\varphi,L_x^{1-\delta_1'}(1+L_x)^{-(\delta_1-\delta_1')}\varphi)
 \end{eqnarray*}
 by Lemmas~\ref{lsub1} and \ref{lsub2} which are applicable if $\delta_1\in[0, 1\wedge n/2\rangle$ and $\delta_1'/\delta_1\in[0,1\rangle$.
 If, however, $n=1$, $\delta_1\in[1/2,1\rangle$ and $\varphi\in C_c^\infty(\Ri^n\backslash\{0\})$  then
 \begin{eqnarray*}
 h_1(\varphi) &\geq&a\,(\varphi',L_{x,D}^{-\delta_1}(1+L_{x,D}^{-\delta_1})^{-1+\delta_1'/\delta_1}\varphi')
 \geq a\,(\varphi,L_{x,N}^{1-\delta_1'}(1+L_{x,N})^{-(\delta_1-\delta_1')}\varphi)
 \end{eqnarray*}
 where the second bound follows from the argument given in Example~6.7 of \cite{ERSZ1}.
 
 The case $\delta_1\leq \delta_1'$ is similar but simpler. 
 It also uses the basic inequality $L_x\geq L_{x,N}$ which then extends to all fractional powers.
 \hfill$\Box$

\bigskip

Next consider the Gru\v{s}in  operator (\ref{ecsg1.2}) on $L_2(\Ri^n\times\Ri^m)$  and let $h$ denote the corresponding quadratic form.
The $\delta$ and $\delta'$  are positive and $\delta_1, \delta_1'\in[0,1\rangle$ but we do not
place any restriction on $\delta_2$ and $\delta_2'$.
Clearly $H\geq a\, H_1(\,=a\,H_1\otimes \one)$  and so the bounds of Proposition~\ref{pcsg3.1} are applicable.
But then one has the following complementary bounds.

\begin{prop}\label{pcsg3.2}
The subelliptic estimate 
$h\geq f$ is valid on $L_2(\Ri^n\times\Ri^m)$
where $f$ is the form of the operator $F$  with
\[
 F\sim L_{x_2}^{\alpha}(1+ L_{x_2})^{\alpha-\alpha'}
\]
where
\[
\alpha=(1-\delta_1)/(1+\delta_2-\delta_1)\;\;\;\;\;\;{ and}\;\;\;\;\;\;
\alpha'=(1-\delta_1')/(1+\delta_2'-\delta_1')\;\;\;.
\]
\end{prop}
\proof\
First one has $h\geq a\, h_\delta$ and after 
 a partial Fourier transformation, i.e., a transformation with respect to the $x_2$ variable, $H_\delta$ 
transforms to an  operator $\widetilde H_\delta$ on $L_2(\Ri^n)$ 
\begin{eqnarray*}
\widetilde H_\delta=H_1+c_{\delta_2, \delta'_2}\,|p_2|^2
\end{eqnarray*}
where $H_1$ is the self-adjoint operator analyzed in Proposition~\ref{pcsg3.1}.
Therefore one may apply the latter proposition to  $H_1$ to bound $\widetilde H_\delta$
below by a differential operator in the $\Ri^n$ variable.
Now to proceed we again use the fractional Hardy inequality but to cover all the relevant cases
we have to pass to a fractional power of $\widetilde H_\delta$.
This is achieved with the aid of  the following simple observation.

\begin{lemma}\label{lsub3}
Let $A$ and $B$ be positive  self-adjoint operators such that the form sum 
$A+B$ is densely defined.
Then the form sum $A^{1/2}+B^{1/2}$ is densely defined and 
\[
(A+B)^{1/2}\geq 2^{-1/2}(A^{1/2}+B^{1/2})
\;\;\;.
\]
Similarly $A^{1/2^n}+B^{1/2^n}$ is densely defined and
\[
(A+B)^{1/2^n}\geq 2^{-1+2^{-n}}(A^{1/2^n}+B^{1/2^n})
\]
for all positive integers $n$.
\end{lemma}
\proof\
First one has $A+B\geq A$.
Hence $(A+B)^{1/2}\geq A^{1/2}$ and $D(A^{1/2})\supseteq D( (A+B)^{1/2})$.
Similarly  $(A+B)^{1/2}\geq B^{1/2}$ and $D(B^{1/2})\supseteq D( (A+B)^{1/2})$.
The first statement then follows because
\[
A+B=(A^{1/2}+B^{1/2})^2/2+(A^{1/2}-B^{1/2})^2/2\geq (A^{1/2}+B^{1/2})^2/2
\;\;\;.
\]
The second statement  follows by iteration.
\hfill$\Box$

\bigskip

The proof of Proposition~\ref{pcsg3.2} now continues by applying Lemma~\ref{lsub3} to deduce that
\[
\widetilde H_\delta^{1/4}
\geq  a\,(H_1^{1/4}
+c_{\delta_2, \delta'_2}^{\;\;\;\;\;1/4}\,|p_2|^{1/2})
\;\;\;.
\]
But if $\delta_1\leq\delta_1'$ and 
 $n\geq2$, or $n=1$ and $\delta_1, \delta_1'\in[0,1/2\rangle$,
then, by Proposition~\ref{pcsg3.1},
\begin{equation}
H_1^{1/4}\geq a\,
L_{x_1}^{\,(1-\delta_1')/4}\,(I+L_{x_1}^{\,(1-\delta_1')/4})^{-1+\sigma}
\label{ecsg3.22}
\end{equation}
where $\sigma=(1-\delta_1)/(1-\delta_1')\in\langle0,1]$.
But $0\leq (1-\delta_1')/4\leq 1/4$.
Therefore 
\[
L_{x_1}^{\,(1-\delta_1')/4}\geq a\, |x_1|^{-(1-\delta_1')/2}
\]
by Lemma~\ref{lsub1}.
Hence
\begin{eqnarray*}
L_{x_1}^{\,(1-\delta_1')/4}\,(I+L_{x_1}^{\,(1-\delta_1')/4})^{-1+\sigma}
&\geq &a\,|{x_1}|^{-(1-\delta_1')/2}(1+|{x_1}|^{-(1-\delta_1')/2})^{-1+\sigma}\\[5pt]
&\geq& a\,|{x_1}|^{-(\delta_1-\delta_1')/2}(1+|{x_1}|)^{-(1-\delta_1)/2}
\end{eqnarray*}
by Lemma~\ref{lsub2}.
Then, however,
\[
\widetilde H_\delta^{1/4}\geq a\,\Big(|{x_1}|^{-(\delta_1-\delta_1')/2}(1+|{x_1}|)^{-(1-\delta_1)/2}+c_{\delta_2, \delta'_2}(x_1)^{1/4}\,|p_2|^{1/2}\Big)
\;\;\;.
\]
But the right hand side is a function of $|{x_1}|$ with  a strictly positive minimum $m$  which 
is estimated by elementary arguments.
The minimum value $m$ is a positive function of $|p_2|$  which then  gives the bound $ \widetilde H_\delta\geq M  I$ with $M=m^4$.

In order to estimate the minimum $m$ we note that the first function in the last estimate, $x_1\mapsto |{x_1}|^{-(\delta_1-\delta_1')/2}(1+|{x_1}|)^{-(1-\delta_1)/2}$,  is decreasing and the second function
in the estimate, $x_1\mapsto c_{\delta_2, \delta'_2}(x_1)^{1/4}\,|p_2|^{1/2}$, is increasing.
Now if $|p_2|$ is small the graphs of the two functions intersect at a unique large value of 
$x_1\sim |p_2|^{-1/(1+\delta_2'-\delta_1')}$.
At this point the value of the sum of the functions is proportional to $|p_2|^{\alpha'/2}$.
Similarly if $|p_2|$ is large the graphs intersect at a small $x_1\sim |p_2|^{-1/(1+\delta_1'-\delta_1)}$
and the minimum value is proportional to $|p_2|^{\alpha/2}$.
Therefore 
\[
m(p_2)\sim |p_2|^{\alpha'/2}(1+|p_2|^2)^{(\alpha-\alpha')/4}\;\;\;\;\;{\rm and }\;\;\;\;\;
M(p_2)\sim |p_2|^{2\alpha'}(1+|p_2|^2)^{\alpha-\alpha'}
\]
which is equivalent to the  bound stated in  the proposition.

If $n=1$ and $\delta_1\in[1/2,1\rangle$ or $\delta_1'\in[1/2,1\rangle$ then the estimation procedure has to be
slightly modified.
Then Proposition~\ref{pcsg3.1}  gives the lower bound
\begin{equation}
H_1^{1/4}\geq a\,
L_{x_1,N}^{\,(1-\delta_1')/4}\,(I+L_{x_1,N}^{\,(1-\delta_1')/4})^{-1+\sigma}
\label{ecsg3.23}
\end{equation}
and $L_{x_1,N}\leq L_{x_1}$.
Nevertheless if $\alpha<1/2$ then $L_{x_1,N}^\alpha\geq a\,L_{x_1}^\alpha$
since the boundary conditions do not affect small fractional powers.
Therefore the bound (\ref{ecsg3.22}) follows from (\ref{ecsg3.23}) by another application of 
Lemma~\ref{lsub2}.

The case $\delta_1\geq\delta_1'$ is similar but simpler.
\hfill$\Box$

\bigskip

Now we are prepared to estimate the crossnorm of the semigroup.

\smallskip

\noindent{\bf Proof of Proposition~\ref{pcsg3.0}}
It follows by definition that $h\sim h_\delta$ and then by  combination of Propositions~\ref{pcsg3.1} and \ref{pcsg3.2} that $h\geq f$ or $h\geq f_N$ where
$f$ is the form of a multiplier $F(L_{x_1},L_{x_2})$ and $f_N$ the form of $F(L_{x_1,N},L_{x_2})$.
Now consider the case $\delta_1\leq\delta_1'$ and assume 
 $n\geq2$, or $n=1$ and $\delta_1, \delta_1'\in[0,1/2\rangle$.
 Then one can apply Lemmas~\ref{lcsg2.1} and \ref{lcsg2.2} with 
 \[
 F(L_{x_1},L_{x_2})=a\,\Big( L_{x_1}^{(1-\delta_1')}(1+ L_{x_1})^{-(\delta_1-\delta_1')}
 +L_{x_2}^{\alpha'}(1+ L_{x_2})^{\alpha-\alpha'}\Big)
 \;\;\;.
 \]
Therefore  $V_F(r)\sim V_{F_1}(r)V_{F_2}(r)$ where 
 \[
 V_{F_1}(r)=|\{p_1: |p_1|^{2(1-\delta_1')}(1+ |p_1|^2)^{-(\delta_1-\delta_1')}<r^2\}|
 \]
 and
 \[
  V_{F_2}(r)=|\{p_2: |p_2|^{2\alpha'}(1+ |p_2|^2)^{\alpha-\alpha'}<r^2\}|
  \;\;\;.
  \]
But  $V_{F_1}(r)\sim r^{n(1/(1-\delta_1'),1/(1-\delta_1))}$.
Similarly $V_{F_2}(r)\sim r^{m(1/\alpha',1/\alpha)}$.
Therefore  $V_F(r)\sim r^{(D',D)}$.
Then the semigroup estimates of Proposition~\ref{pcsg3.0} follow from Lemma~\ref{lcsg2.2}.

The argument is similar if $\delta_1\geq\delta_1'$ but one uses the second estimate of Proposition~\ref{pcsg3.1}.
Finally if $n=1$ and $\delta_1\in[1/2,1\rangle$ or $\delta_1'\in[1/2,1\rangle$ one can make an identical argument
using the last statement of Proposition~\ref{pcsg3.1} and Lemma~\ref{lcsg2.3} which deals with the Neumann multipliers.
\hfill$\Box$

\bigskip

\begin{remarkn}
The arguments we have given for subellipticity estimates on $L_2(\Ri^d\,;dx)$ also extend to weighted
spaces such as  $L_2(\Ri^d\,;|x|^\beta dx)$.
In this extension the Hardy inequality is replaced by the  Caffarelli--Kohn--Nirenberg inequalities \cite{CKN}.
\end{remarkn}

\section{Comparison of kernels}\label{Scsg4}

Proposition~\ref{pcsg3.0} gives  bounds on the semigroup kernel associated with the Gru\v{s}in  operator  which are uniform over $\Ri^n\times \Ri^m$.
In this section we  develop a method for transforming these bounds into bounds which better reflect the
spatial behaviour,  bounds expressed in terms of the corresponding Riemannian geometry.
In particular we establish a comparison between the Gru\v{s}in kernel and the kernel of a closely related non-degenerate operator.
Our arguments are based on wave equation techniques and the initial  problem is to establish that the wave
equation has a finite speed of propagation when measured with respect to the Riemannian distance
(\ref{ecsg1.21}).
We discuss this problem in the general context of the elliptic operators defined by the relaxation $h_0$ of the form
$h$ corresponding to an elliptic operator $H$ given by (\ref{ecsg1.1}).

First recall that the Riemannian distance is given  by
\begin{equation}
d(x\,;y)=\sup_{\psi\in D}\,(\psi(x)-\psi(y))
\label{ecsg4.00}
\end{equation}
for all $x,y\in\Ri^d$
where 
 \[
D=\{\psi\in W^{1,\infty}(\Ri^d):\sum^d_{i,j=1}c_{ij}\,(\partial_i\psi)(\partial_j\psi)\leq1\}
\;\;\;.
\]
Secondly, introduce the corresponding set-theoretic distance by
\begin{equation}
d(A\,;B)=\inf_{x\in A,\,y\in B}d(x\,;y)
\label{ecsg4.100}
\end{equation}
where $A$ and $B$ are general measurable sets.

Our first aim is to establish the following basic estimates.

\begin{prop}\label{pscsg4.0}
Let $S^{(0)}$ denote the semigroup generated by the self adjoint operator $H_0$ associated with the 
relaxation $h_0$ of the elliptic form $h$ with coefficients $C=(c_{ij})$.
Then for each pair of open subsets $A,B$ of $\Ri^d$
\begin{equation}
|(\varphi_A, S^{(0)}_t\varphi_B)|\leq e^{-d(A;B)^2(4t)^{-1}}\|\varphi_A\|_2\|\varphi_B\|_2
\label{epadty3.0}
\end{equation}
for all $\varphi_A\in L_2(A)$, $\varphi_B\in L_2(B)$ and $t>0$ with the convention
$e^{-\infty}=0$.
Moreover, the corresponding wave equation has a finite speed of propagation in the sense that
\begin{equation}
(\varphi_A,{\cos(tH_0^{1/2})}\varphi_B)=0
\label{epadty3.1}
\end{equation}
for all $\varphi_A\in L_2(A)$, $\varphi_B\in L_2(B)$ and all $t$ with  $|t|\leq d(A\,;B)$.
\end{prop}
\proof\
First let $\psi\in D$ and introduce the one-parameter family of multiplication operators
$\rho\to U_\rho=e^{\rho\psi}$ on $L_2(\Ri^d)$.
Then $\|U_\rho S^{(0)}_tU_\rho^{-1}\|_{2\to2}\leq e^{\rho^2t}$.
This is a standard estimate for strongly elliptic operators which extends to general elliptic operators
(see, for example, the proof of Proposition~3.1 in \cite{ERSZ1}).
Now
\begin{eqnarray*}
|(\varphi_A,S^{(0)}_t\varphi_B)|&=&
|(U_\rho^{-1}\varphi_A,(U_\rho S^{(0)}_tU_\rho^{-1})U_\rho\varphi_B)|\\[5pt]
&\leq&e^{\rho^2t}\,\|U_\rho^{-1}\varphi_A\|_2\,\|U_\rho\varphi_B\|_2
\leq e^{-\rho d_\psi(A\,;B)}e^{\rho^2t}\,\|\varphi_A\|_2\,\|\varphi_B\|_2
\end{eqnarray*}
where
\[
d_\psi(A\,;B)=\inf_{x\in A,\,y\in B}(\psi(x)-\psi(y))
\;\;\;.
\]
Therefore, optimizing over $\psi$  and $\rho$  one has
\begin{equation}
|(\varphi_A,S^{(0)}_t\varphi_B)|\leq  e^{-\hat d(A;B)^2(4t)^{-1}}\|\varphi_A\|_2\|\varphi_B\|_2
\label{epadty3.9}
\end{equation}
where 
\begin{equation}
\hat d(A\,;B)=\sup_{\psi\in D}d_\psi(A\,;B)
\;\;\;.
\label{epadty3.10}
\end{equation}
These estimates are valid for all measurable $A, B$ and all $\varphi_A\in L_2(A),\varphi_B\in L_2(B)$.

Since $d_\psi(A\,;B)\leq \psi(x)-\psi(y)$ for all $x\in A$ and $y\in B$ it follows that 
$\hat d(A\,;B)\leq d(A\,;B)$ again for all measurable $A$ and $B$.
But the latter inequality has a partial converse.

\begin{lemma}\label{lcsg4.0}
If  $A$ and $B$ are compact subsets then $\hat d(A\,;B)= d(A\,;B)$ .
\end{lemma}
\begin{remarkns} 1. The proof is an interplay between compactness and continuity. 
It uses very little structure of the set $D$.
Indeed it suffices for the proof that $\varphi\in D$ and $c\in \Ri$ imply $-\varphi+c\in D$
and $\varphi_1,\varphi_2\in D$ imply $\varphi_1\vee \varphi_2, \varphi_1\wedge\varphi_2 \in D$.
These properties are easily verified.

\hspace{1.2cm}2. It is not necessarily the case that $\{\varphi_n\}_{n\geq1}\in D$ implies $\sup_{n\geq1} \varphi_n\in D$ or $\inf_{n\geq1}\varphi_n\in D$.
Proposition~6.5 in \cite{ERSZ1},  and its  proof, give an example of a  decreasing sequence $\chi_n$ of functions in $D$ such that $\inf_n\chi_n\not\in W^{1,\infty}$.
\end{remarkns}
\noindent{\bf Proof of Lemma~\ref{lcsg4.0}}
Since $\hat d(A\,;B)\leq d(A\,;B)$  for all measurable $A$ and $B$
it suffices to prove  $\hat d(A\,;B)\geq d(A\,;B)$ for $A$ and $B$ compact.

Fix $x\in A$ and $y\in B$.
Then for each $\varepsilon>0$ there is a $\psi_{x,y}\in D$ such that $\psi_{x,y}(x)=0$ and  
$\psi_{x,y}(y)\geq d(x\,;y)-\varepsilon/4$.
Next since $\psi_{x,y}$ is continuous there exists an open neighbourhood  $U_y$ of $y$ such that
\[
\psi_{x,y}(z)\geq d(x\,;y)-\varepsilon/2\geq d(A\,;B)-\varepsilon/2
\]
for and $z\in U_y$.
Then $B\subset \bigcup_{y\in B}U_y$ and since $B$ is compact there exist
$y_1,\ldots,y_n\in B$ such that $B\subset \bigcup_{k=1}^nU_{y_k}$.
Set $\psi_x=\sup_{1\leq k\leq n}\psi_{x,y_k}$ then $\psi_x\in D$, $\psi_x(x)=0$ and
\[
\inf_{z\in B}\psi_x(z)\geq \min_{1\leq k\leq n}d(x\,;y_k)-\varepsilon/2\geq d(A\,;B)-\varepsilon/2
\;\;\;.
\]
But since $\psi_x$ is continuous there is an open neighbourhood $U_x$ of $x$ such that $\psi_x(z)\leq \varepsilon/2$ for all $z\in U_x$.
Then by repeating the above covering  argument one can select  $x_1,\ldots,x_m\in A$ 
such that $\psi=\inf_{1\leq k\leq m}\psi_{x_k}$ satisfies $\psi(z)\leq \varepsilon/2$ for all $z\in A$.
In addition one still has $\psi(z)\geq d(A\,;B)-\varepsilon/2$ for all $z\in B$.
Therefore
\[
\inf_{x\in A,y\in B}(\psi(y)-\psi(x))\geq d(A\,;B)-\varepsilon
\;\;\;.
\]
Thus  $\hat d(A\,;B)\geq d(A\,;B)-\varepsilon$.
Since $\varepsilon>0$ is arbitrary one deduces that  $\hat d(A\,;B)\geq d(A\,;B)$ for $A$ and $B$ compact. \hfill$\Box$

\bigskip

Now we can complete the proof of Proposition~\ref{pscsg4.0}

\smallskip

\noindent{\bf End of  proof of Proposition~\ref{pscsg4.0}}
Let  $\varphi_1\in L_2(A)$ have compact support $U_1\subset A$ and $\varphi_2\in L_2(B)$
have compact support $U_2\subset B$.
Then it follows from (\ref{epadty3.9}) and Lemma~\ref{lcsg4.0} that 
\begin{eqnarray*}
|(\varphi_1,S^{(0)}_t\varphi_2)|&\leq & 
e^{-\hat d(U_1;U_2)^2(4t)^{-1}}\|\varphi_1\|_2\|\varphi_2\|_2\\[5pt]
&=&e^{- d(U_1;U_2)^2(4t)^{-1}}\|\varphi_1\|_2\|\varphi_2\|_2
\leq e^{- d(A;B)^2(4t)^{-1}}\|\varphi_1\|_2\|\varphi_2\|_2
\;\;\;.
\end{eqnarray*}
But if $A$ is  open  then the functions of compact support in $L_2(A)$ are dense
and similarly for $B$.
Therefore the first statement  (\ref{epadty3.0}) of the proposition follows by continuity.
The second statement (\ref{epadty3.1})  is a direct consequence (see, for example, \cite{ERSZ1} Lemma~3.3).
\hfill$\Box$

\bigskip

An equivalent way of expressing the finite speed of  propagation  is the following.
\begin{lemma} \label{lpre2.11}
Let $A$  be an  open subset and $F$ a closed subset   with $A\subset F$.
Then 
\[
{\cos(tH_0^{1/2})}L_2(A)\subseteq L_2(F)
\]
for all $t\in\Ri$ with $|t|\leq d(A\,;F^{\rm c})$.
\end{lemma}
\proof\
Let $B$ be an  open subset of $ F^{\rm c}$. Then $d(A\,;F^{\rm c})\leq d(A\,;B)$.
Therefore $L_2(B)\perp {\cos(tH_0^{1/2})}L_2(A)$ for all $t\in\Ri$ with $|t|\leq d(A\,;F^{\rm c})$
by (\ref{epadty3.1}).
Hence $L_2(F^{\rm c})\perp {\cos(tH_0^{1/2})}L_2(A)$ and one must have 
${\cos(tH_0^{1/2})}L_2(A)\subseteq L_2(F)$.\hfill$\Box$

\bigskip

  The subsequent comparison theorem depends on a generalization of the propagation property which
  emphasizes the local nature.
  As a preliminary let $h_1$ and $h_2$ be two  elliptic forms with coefficients $C_1=(c^{(1)}_{ij})$ and $C_2=(c^{(2)}_{ij})$ and let $H_{1,0}$ and $H_{2,0}$ denote the corresponding relaxations.
  Moreover, assume that $C_1\geq C_2$ and $C_1\geq\mu I>0$.
  In particular $D(h_1)\subseteq  D(h_2)\subseteq W^{1,2}(\Ri^d)$ and $d_1(x\,;y)\leq d_2(x\,;y)$ for all $x,y$ where $d_1$ and $d_2$ denote the Riemannian distances associated with $C_1$ and $C_2$.
  Let $U=\supp (C_1-C_2)$ and set $F=\overline{U^{\rm c}}$.
  
  \begin{lemma}\label{lcsg4.2} If $A$ is an open subset of the closed subset $F$ then
  \begin{equation}
\cos(tH_{1,0}^{1/2})\varphi_A=\cos(tH_{2,0}^{1/2})\varphi_A
\label{ecsg4.5}
\end{equation}
for all $\varphi_A\in L_2(A)$ and all $t\in \Ri$ with $|t|\leq d_1(A\,;U)$.
\end{lemma}
\proof\
If $C_1$ and $C_2$ are strongly elliptic this result follows from the proof  of Proposition~3.15 in 
\cite{ERS4}.
The extension to the more general situation can then be made by approximation as follows.

Let $N>\varepsilon>0$.
Set $C_{1,N,\varepsilon}=(C_1\wedge N I)+\varepsilon I$ and $C_{2,N,\varepsilon}=(C_2\wedge N I)+\varepsilon I$.
Then the lemma is valid for the corresponding strongly elliptic operators $H_{1,N,\varepsilon}$,
$H_{2,N,\varepsilon}$.
But as $N\to\infty$ these operators converge in the strong resolvent sense to 
$H_{1,\varepsilon}=H_1+\varepsilon I$ and $H_{2,\varepsilon}=H_2+\varepsilon I$, respectively.
Since $H_{1,\varepsilon}\geq H_{1,N,\varepsilon}\geq H_{2,N,\varepsilon}$  it follows that 
\[
\cos(tH_{1,\varepsilon}^{1/2})\varphi_A=\cos(tH_{2,\varepsilon}^{1/2})\varphi_A
\]
for all $t\in\Ri$ with $|t|\leq d_{1,\varepsilon}(A\,;U)$ where $d_{1,\varepsilon}$
denotes the Riemannian distance associated with $C_{1,\varepsilon}$.

Finally it follows that $H_{1,\varepsilon}$ and $H_{2,\varepsilon}$ converge in the strong resolvent sense to $H_{1,0}$ and $H_{2,0}$, respectively, as $\varepsilon\to0$.
Moreover, since $C_1\leq C_{1,\varepsilon}\leq (1+\varepsilon\mu^{-1}) \,C_1$ it follows that
$(1+\varepsilon\mu^{-1})^{-1/2}d_1(A\,;B)\leq d_{1,\varepsilon}(A\,;B)\leq d_1(A\,;B)$
for all measurable $A$, $B$.
Hence $d_{1,\varepsilon}(A\,;B)\to d_1(A\,;B)$ as $\varepsilon\to0$.
Therefore the statement of the lemma follows in the limit.
\hfill$\Box$

  \bigskip

Now we are prepared to establish the principal comparison result.
In the sequel $S^{(1,0)}$, $S^{(2,0)}$ denote the semigroups generated by 
$H_{1,0}$, $H_{2,0}$  and  $K^{(1,0)}$, $K^{(2,0)}$ denote the corresponding
kernels.

\begin{thm}\label{tcsg3.1}
Adopt the foregoing notation and assumptions.
Let 
$\chi_A$  denote the characteristic function of the open subset 
$A$ of $ F$.
Set 
\[
M_N(t)=\|\chi_A(I+t^2\,H_{1,0})^{-N}\chi_A\|_{1\to\infty}+\|\chi_A(I+t^2\,H_{2,0})^{-N}\chi_A\|_{1\to\infty}
\]
for $N\in\Ni$.

Then there is an $a_N>0$  such that 
\begin{equation}
\sup_{x,y \in A} |K^{(1,0)}_t(x\,;y)- K^{(2,0)}_t(x\,;y)| \leq a_N\,
M_N(t/\rho )\,(\rho^2/t)^{-1/2}
 \,e^{-\rho^2/(4t)}
\label{ecsg3.2}
\end{equation}
for all $t>0$ where $\rho =d_1(A\,;U)$.
\end{thm}

Note that there is no reason that $M_N$ is finite.
Subsequently we give conditions which
ensure that $M_N$ is indeed finite for large $N$.

The proof of the theorem is based on the estimates developed in \cite{Sik} and \cite{Sik3}.
As a preliminary we need some properties of functions of the operators $H_{0, 1}$ and 
$H_{2,0}$.

First, let $\Psi$ be an  even bounded Borel function with
Fourier transform 
$\widetilde \Psi$ satisfying $\supp \widetilde\Psi\subseteq [-\rho, \rho]$ where $\rho>0$.
Then for each pair of  open subsets $B$ and $C$ 
\begin{equation}
(\varphi_B,\Psi(H_{1,0}^{1/2})\varphi_C)=0=(\varphi_B,\Psi(H_{2,0}^{1/2})\varphi_C)
\label{ecsg3.4}
\end{equation}
for all $\varphi_B\in L_2(B)$, $\varphi_C\in L_2(C)$ where  
$\rho\leq d_1(B\,;C)$.
This follows for $H_{1,0}$ from the representation
\[
 \Psi(H_{1,0}^{1/2})
=(2\pi)^{-1/2}\int_{\Ri}dt\, \widetilde \Psi(t) \exp(itH_{1,0}^{1/2})
=(2\pi)^{-1/2}\int_{-\rho}^\rho dt\, \widetilde \Psi(t) \cos(tH_{1,0}^{1/2})\;\;\;.
\]
and   condition (\ref{epadty3.1}).
The argument for $H_{2,0}$ is similar.

Secondly, we have the key lemma.
\begin{lemma}\label{lcsg3.1}
Let $\Psi$ be an  even bounded Borel function with
Fourier transform 
$\widetilde \Psi$ satisfying $\supp \widetilde\Psi\subseteq [-2\rho ,2\rho ]$.
Then
\[
\chi_A \Psi(H_{1,0}^{1/2})\chi_A=\chi_A \Psi(H_{2,0}^{1/2})\chi_A
\;\;\;.
\]
\end{lemma}
\proof\
Since  $A\subset F$
and  $\rho =d_1(A\,;U)$ one has 
\begin{equation}
\cos(tH_{1,0}^{1/2})\varphi_A=\cos(tH_{2,0}^{1/2})\varphi_A
\label{ecsg3.5}
\end{equation}
for all  $\varphi_A\in L_2(A)$ and $t\in \Ri$ with $|t|\leq \rho$ by Lemma~\ref{lcsg4.2}.

Next remark that 
\begin{eqnarray*}
(\psi_A,\cos(2tH_{1,0}^{1/2})\varphi_A)-(\psi_A,\cos(2tH_{2,0}^{1/2})\varphi_A)\\[5pt]
&&\hspace{-5cm}=
2\,\Big((\cos(tH_{1,0}^{1/2})\psi_A,\cos(tH_{1,0}^{1/2})\varphi_A)-
(\cos(tH_{2,0}^{1/2})\psi_A,\cos(tH_{2,0}^{1/2})\varphi_A)\Big)=0
\end{eqnarray*}
for $|t|\leq \rho$ and $\psi_A,\varphi_A\in L_2(A)$.
Then, however, one has
\begin{eqnarray*}
(\psi, \chi_A \Psi(H_{1,0}^{1/2})\chi_A\varphi)
&=&(2\pi)^{-1/2}\int_{-2\rho}^{2\rho}dt\, \widetilde \Psi(t) (\chi_A\psi, \cos(tH_{1,0}^{1/2})\chi_A\varphi)\\[5pt]
&=&(2\pi)^{-1/2}\int_{-2\rho}^{2\rho}dt\, \widetilde \Psi(t) (\chi_A\psi,\cos(tH_{2,0}^{1/2})\chi_A\varphi)=
(\psi, \chi_A \Psi(H_{2,0}^{1/2})\chi_A\varphi)
\end{eqnarray*}
for all $\varphi\in L_2(\Ri^d)$.
\hfill$\Box$

\bigskip

Now we are prepared to prove the theorem.
We use the notation $K_T$ for the kernel of an operator $T$.

\smallskip

 \noindent{\bf Proof of Theorem~\ref{tcsg3.1}}
 Let  $\psi \in C^{\infty}(\Ri)$ be an increasing function with
\[
\psi(u) = \left\{ \begin{array}{ll}
    0 & \mbox{ if $u \le -1$}\\
    1  & \mbox{ if $u \ge -1/2$}\; .
           \end{array}
    \right.
\]
Then for $s>1$ define the family of functions $\varphi_s$ such that 
\[
\varphi_s(u)= \psi(s(|u|-s))\;\;\;.
\]
Next  define functions $\widetilde\Phi_s$ and $\widetilde\Psi_s$ by 
\begin{equation}
\widetilde\Phi_s(u)
= (4\pi)^{-1/2}\exp {({-u^2/4})}-\widetilde\Psi_s(u)
= \varphi_s(u) \,(4\pi)^{-1/2}\exp {({-u^2/4})}
\;\;\;.
\end{equation}
Then  the inverse  Fourier transforms satisfy 
${\Phi_s}(\lambda)+{\Psi_s}(\lambda)=(2\pi)^{-1/2}\,\exp(-\lambda^2)$
and
\begin{equation}
S^{(i,0)}_t=\exp(-tH_{i,0})={\Phi_s}((tH_{i,0})^{1/2})+{\Psi_s}((tH_{i,0})^{1/2})
\label{ggg}
\end{equation}
for $i=1,2$.
Integration by parts $2N$ times yields
\begin{eqnarray*}
\int du\,e^{-u^2/4}e^{-iu\lambda}\,\varphi_s(u)=
\int du\,e^{-u^2/4-i\lambda u}\,\underbrace{\Big(\frac{1}{u/2+i\lambda}%
\Big(\ldots\Big(\frac{1}{u/2+i\lambda}\varphi_s(u)\Big)^{\prime}\ldots\Big)'\Big)%
^{\prime}}_{2N}\;\;\;.
\end{eqnarray*}
Hence for any $N\in\Ni$ and $s>1$ there is an $a_N>0$ such that
\begin{equation}\label{osz}
|{\Phi_s}(\lambda)| \le
a_N\,\frac{1}{s\,(1+\lambda^2/s^2)^{N}}\,e^{-s^2/4} \;\;\;,\label{osz1}
\end{equation}
with the value of $a_N$ depending only on $N$.

Next note that
supp~$\widetilde\Psi_s \subseteq [-s+(2s)^{-1},s-(2s)^{-1}]\subset [-s,s]$.
So setting $s_{\rho}= 2\,\rho\,t^{-1/2} $ one has 
\[
\chi_A{{\Psi_{s_{\rho}}}((tH_{1,0})^{1/2})}\chi_A=
\chi_A{{\Psi_{s_{\rho}}}((tH_{2,0})^{1/2})}\chi_A
\]
by Lemma~\ref{lcsg3.1} and rescaling with $t^{1/2}$.
Hence  one deduces from (\ref{ggg}) that 
\begin{eqnarray}
\chi_A(S^{(1,0)}_t-S^{(2,0)}_t)\chi_A&=&
\chi_A{{\Phi_{s_{\rho}}}((tH_{1,0})^{1/2}))}\chi_A-\chi_A{{\Phi_{s_{\rho}}}((tH_{2,0})^{1/2}))}\chi_A
\;\;\;.
\label{incl}
\end{eqnarray}
Therefore
\begin{eqnarray*}
\sup_{x,y \in A} |K^{(1,0)}_t(x\,;y)- K^{(2,0)}_t(x\,;y)|&&\\[5pt]
& &\hspace{-2cm}{}\leq
\|\chi_A{\Phi_{s_{\rho}}}((tH_{1,0})^{1/2})\chi_A\|_{1\to\infty}+
\|\chi_A{\Phi_{s_{\rho}}}((tH_{2,0})^{1/2})\chi_A\|_{1\to\infty}
\;\;\;.
\end{eqnarray*}

Now let $\Theta_{s_{\rho}}$ be a possibly complex function such that
$\Theta_{s_{\rho}}(\lambda)^2={\Phi_{s_{\rho}}}(t^{1/2}\lambda)$.
Then
\begin{eqnarray*}
\|\chi_A{\Phi_{s_{\rho}}}((tH_{i,0})^{1/2}))\chi_A\|_{1\to\infty}&=&
(\|\chi_A\,\Theta_{s_{\rho}}((H_{i,0})^{1/2})\|_{2\to\infty})^2\\[15pt]
&&\hspace{-3cm}\leq (\|(I+t^2H_{i,0}/\rho^2)^{N/2}\,\Theta_{s_{\rho}}((H_{i,0})^{1/2})\|_{2\to2})^2\,(\|\chi_A(I+t^2H_{i,0}/\rho^2)^{-N/2}\|_{2\to\infty})^2
\;\;\;.
\end{eqnarray*}
But
\begin{eqnarray*}
\left|\Theta_{s_{\rho}}(\lambda)
(1+t^2\lambda^2/\rho^2)^{N/2}\right|^2
&=&
\left|{\Phi_{s_{\rho}}}(t^{1/2}\lambda)
(1+t(t^{1/2}\lambda)^2/\rho^2)^{N}\right|\\[5pt]
&=&
\left|{\Phi_{s_{\rho}}}(t^{1/2}\lambda)
(1+4(t^{1/2}\lambda)^2/{s_\rho}^2)^{N}\right|\;\;\;.
\end{eqnarray*}
Therefore
\begin{eqnarray*}
(\|(I+t^2H_{i,0}/\rho^2)^{N/2}\,\Theta_{s_{\rho}}((H_{i,0})^{1/2})\|_{2\to2})^2&=&
\sup_{\lambda \ge 0}\left|\Theta_{s_{\rho}}(\lambda)
(1+t^2\lambda^2/\rho^2)^{N/2}\right|^2\\[5pt]
&=&\sup_{\lambda \ge 0}
\left|{\Phi_{s_{\rho}}}(t^{1/2}\lambda)
(1+4(t^{1/2}\lambda)^2/{s_\rho}^2)^{N}\right|\\[5pt]
&\leq& a_N\,(s_\rho/2)^{-1}\,e^{-s_\rho^2/16}=a_N\,(\rho^2/t)^{-1/2}\,e^{-\rho^2/(4t)}
\;\;\;.
\end{eqnarray*}

Combining these estimates gives
\begin{eqnarray*}
\sup_{x,y \in A} |K^{(1,0)}_t(x\,;y)- K^{(2,0)}_t(x\,;y)|
&=&a_N\,M_N(t/\rho)\,(\rho^2/t)^{-1/2}\,e^{-\rho^2/(4t)}
\end{eqnarray*}
which establishes the statement of the theorem.
\hfill$\Box$

\bigskip

The statement of the theorem can be reformulated in terms of {\it a priori}
bounds on the semigroups $S^{(1,0)}$ and $S^{(2,0)}$.
But this requires a uniform estimate on the crossnorms $\|S^{(i,0)}\|_{1\to\infty}$.

 If $h$ is strongly elliptic then  $\|S^{(0)}_t\|_{1\to\infty}\leq a\,t^{-d/2}$ for all $t>0$ and we will
 assume  analogous bounds 
 \[
 \|S^{(0)}_t\|_{1\to\infty}\leq a\,V(t)^{-1}
 \]
 where $V$ is a positive increasing function which satisfies the doubling property
 \begin{equation}
V(2t)\leq a\,V(t)
\label{epre1.1}
\end{equation}
for some $a>0$ and all $t>0$.
It  follows from (\ref{epre1.1})  that there are $a, \widetilde D>0$  such that 
\begin{equation}
V(s)\leq a \,(s/t)^{\widetilde D}\,V(t)
\label{epre1.17}
\end{equation}
for all $s\geq t>0$.
The parameter $\widetilde  D$ is the doubling dimension, although it need not be an integer.

\begin{cor}\label{ccsg3.1}
Adopt the hypotheses and notation of Theorem~$\ref{tcsg3.1}$.
Let $V$ be a positive increasing function which satisfies the doubling property~$(\ref{epre1.1})$.
Assume that 
\[
\|\chi_A\,S^{(1,0)}_t\,\chi_A\|_{1\to\infty}\vee \|\chi_A\,S^{(2,0)}_t\,\chi_A\|_{1\to\infty}\leq V(t)^{-1}
\]
for all $t>0$.

Then there is an $a>0$ such that 
\begin{equation}
\sup_{x,y \in A} |K^{(1,0)}_t(x\,;y)- K^{(2,0)}_t(x\,;y)| \leq a\,
V(t^2/\rho^2)^{-1}\,(\rho^2/t)^{-1/2}
 \,e^{-\rho^2/(4t)}
\label{ecsg3.20}
\end{equation}
for all $t>0$, where $\rho=d_1(A\,;U)$.
\end{cor}
\proof\
It follows by  Laplace transformation that
\begin{eqnarray*}
\|\chi_A (I+t^2\,H_0)^{-N}\chi_A\|_{1\to\infty}&\leq& 
{{1}\over{(N-1)!}}\int^\infty_0ds\,s^{N-1}e^{-s}\|\chi_A\,
 S^{(0)}_{st^2}\,\chi_A\|_{1\to\infty}\\[5pt]
&\leq &{{1}\over{(N-1)!}}\int^\infty_0ds\,s^{N-1}e^{-s}\,V(st^2)^{-1}\leq a_N\,V(t^2)^{-1}
\end{eqnarray*}
for all $t>0$ with $a_N$ finite if $N>\widetilde  D$.
The last step uses the doubling property (\ref{epre1.17}) in the form $V(t^2)\leq \rho\,s^{-\widetilde D}\,V(st^2)$ for $s\leq 1$.
Therefore $M_N(t)\leq 2\,a_N\,V(t^2)^{-1}$ and the statement of the corollary is an immediate
consequence of Theorem~\ref{tcsg3.1}.
\hfill$\Box$

\bigskip

In Section~\ref{Scsg6} we apply Corollary~\ref{ccsg3.1} to the Gru\v{s}in operator $H$.
Then we set $H_2=H$ and $H_1$ a Gru\v{s}in operator  with no local degeneracies
but with the same growth properties for $|x_1|\geq 1$.

\section{Volume estimates}\label{Scsg5}

In this section we return to the analysis of the general Gru\v{s}in
operator $H$.
Our aim is to calculate the Riemannian distance  $d(\cdot\,;\cdot)$, given by
 (\ref{ecsg1.21}),
 and   the  volume of the  corresponding balls
$B(x\,;r)=\{y\in\Ri^n\times\Ri^m: d(x\,;y)<r\}$, up to equivalence.
 In Section~\ref{Scsg3} we established uniform bounds on the crossnorm
$\|S_t\|_{1\to\infty}$
 on the semigroup generated by the relaxation of $H$ and these
automatically give bounds
 on the semigroup kernel~$K_t$ which are uniform over $\Ri^n\times \Ri^m$.
 But in Section~\ref{Scsg6} we will improve these uniform bounds, with the
help of the comparison
 results of Section~\ref{Scsg4} to obtain bounds
$K_t(x\,;y)\leq
a\,(|B(x;t^{1/2})|\,|B(y;t^{1/2})|)^{-1/2}$
in terms of the volume of the balls.

First, remark that if $C_1$, $C_2$ are two positive symmetric matrices whose entries
are measurable functions and $d_1(\cdot\,;\cdot)$, $d_2(\cdot\,;\cdot)$ the corresponding distances
then $C_1\sim C_2$ implies  $d_1(\cdot\,;\cdot)\sim d_2(\cdot\,;\cdot)$.
In particular if $a\,C_1\leq C_2\leq b\,C_1$ then $b^{-1/2}d_1(x\,;y)\leq d_2(x\,;y)\leq a^{-1/2}d_1(x\,;y)$
for all $x,y\in \Ri^d$.
Moreover, the corresponding balls $B_1$, $B_2$ satisfy
\[
B_2(x\,;b^{-1/2}r)\subseteq B_1(x\,;r)\subseteq B_2(x\,;a^{-1/2}r)
\]
for all $x\in \Ri^d$ and all $r>0$.
Then $|B_1|\sim |B_2|$ and if $|B_1|$ satisfies the doubling property with doubling  dimension $\widetilde D$ then $|B_2|$
also satisfies the property with the same dimension, and conversely.

Secondly, in considering the Gru\v{s}in operator the coefficient matrix $C\sim C_\delta$
and to calculate the Riemannian distance, up to equivalence, we may make a convenient choice of the
$C_\delta$.
In particular we may choose $C_\delta$ such that its entries $c_{\delta_1,\delta'_1}$,  $c_{\delta_2,\delta'_2}$ are continuous functions over $\Ri^n$.
In fact we may assume $c_{\delta_1,\delta'_1}\in C^{\delta_1}(\Ri^n)$ and 
 $c_{\delta_2,\delta'_2}\in C^{\delta_2}(\Ri^n)$.
 The continuity of the coefficients then allow us to appeal to path arguments
 in the computation of the Riemannian distance.
 Specifically the Riemannian distance defined by (\ref{ecsg1.21}) is equivalent to the shortest distance
 of paths measured with respect to a continuous choice of $C_\delta$.

Thirdly, recall that  we always assume that
$0 \leq \delta_1,\delta_1'<1$ and $0 \leq \delta_2,\delta'_2$.
Next   let  $B(x_1,x_2\,;r)$  denote the ball with centre $x=(x_1,x_2)$
and    $|B(x_1,x_2\,;r)|$  its volume.
Further let  $D=(n+m(1+\delta_2-\delta_1))(1-\delta_1)^{-1}$
and  $D'=(n+m(1+\delta'_2-\delta_1'))(1-\delta_1')^{-1}$
denote the parameters occurring in the uniform bounds of
Proposition~\ref{pcsg3.0}.
Next define the function   $\Delta_\delta$ by the formula
\begin{equation}\label{csgdelta}
  \Delta_\delta (x_1,x_2\,;y_1,y_2)=   \left\{ \begin{array}{llll}
{|x_2-y_2|}/{(|x_1|+|y_1|)^{(\delta_2,\delta_2')}}
 & \mbox{ if $|x_2-y_2| \leq  (|x_1|+|y_1|)^{(\rho,\rho')}  $}\\[8pt]
  {|x_2-y_2|^{(1-\gamma,1-\gamma')}  }
&     \mbox{ if
$ |x_2-y_2|\geq (|x_1|+|y_1|) ^{(\rho,\rho')}   $}  
           \end{array}
    \right.
\end{equation}
where $\rho=1+\delta_2-\delta_1$,  $\rho'=1+\delta_2'-\delta_1'$,
$\gamma =\delta_2/\rho$ and 
$\gamma' =\delta'_2/\rho'$.
Now  set
$$
  D_\delta (x_1,x_2\,;y_1,y_2)=
{|x_1-y_1|}/{(|x_1|+|y_1|)^{(\delta_1,\delta_1')}}
+\Delta_\delta (x_1,x_2\,;y_1,y_2).
$$
Note that
$$
 \Delta_\delta (x_1,x_2\,;y_1,y_2)\sim
\frac{|x_2-y_2|}{ (|x_1|+|y_1|)^{(\delta_2,\delta_2')}+
|x_2-y_2|^{(\gamma,\gamma')} }.
$$
Note also that if  $ (|x_1|+|y_1|) ^{(\rho,\rho')}= |x_2-y_2|   $
then
$$
\frac{|x_2-y_2|^{\phantom{(\delta_2,\delta_2')}}}{(|x_1|+|y_1|)^{(\delta_2,\delta_2')}}
=\frac{|x_2-y_2|^{\phantom{(\gamma,\gamma')}}}{|x_2-y_2|^{(\gamma,\gamma')}
}=|x_2-y_2|^{(1-\gamma,1-\gamma')}
$$
so $ \Delta_\delta$ is a continuous function of  the variables
$x_1,x_2,y_1,y_2$.
The main result  in this section is the  following.

\begin{prop}\label{pcsg5.1}
Consider the general Gru\v{s}in operator with coefficients $C\sim C_\delta$.
If $d_\delta$ is the Riemannian distance defined by
$(\ref{ecsg1.21})$ with the coefficients $C_\delta$  then
\[
d_\delta(x_1,x_2\,;y_1,y_2)\sim D_\delta (x_1,x_2\,;y_1,y_2)
\;\;\;.
\]
Moreover, the volume of the corresponding  balls satisfy
\begin{equation}\label{evol}
|B(x_1,x_2;r)|\sim \left\{ \begin{array}{llll}
 r^{(D,D')} & \mbox{ if $r\geq |x_1|^{(1-\delta_1,1-\delta_1')} $}\\[5pt]
   r^{n+m}|x_1|^{(\beta,\beta')}
& \mbox{ if $ r \leq |x_1|^{(1-\delta_1,1-\delta_1')} $}
           \end{array}
    \right.
\end{equation}
where
$\beta=n\delta_1+m\delta_2$ and $\beta'=n\delta'_1+m\delta_2'$.
\end{prop}
\proof\
First note that the coefficients of  $C_\delta$ do not depend on $x_2$. 
Hence
\[
d_\delta(x_1,x_2\,;y_1,y_2)=d_\delta(x_1,0\,;y_1,y_2-x_2)
\;\;\;.
\]
Without lost of generality one may  assume  $|x_1| \leq |y_1|$
and so $|y_1|\sim |x_1|+|y_1|$. 
We adopt  this convention throughout the remainder of the  proof.
Next by the  triangle inequality
\begin{equation}
d_\delta(x_1,0\,;y_1,y_2-x_2)\leq
d_\delta(x_1,0\,;y_1,0)+ d_\delta(y_1,0\,;y_1,y_2-x_2)
\;\;\;.\label{egr5.1}
\end{equation}
Now we  argue that the first term on the right hand side of (\ref{egr5.1})  satisfies the estimate
\[
d_\delta(x_1,0\,;y_1,0)\leq a\,
\frac{|x_1-y_1|^{\phantom{(\delta_1,\delta_1')}}}
{(|x_1|+|y_1|)^{(\delta_1,\delta_1')}}
\;\;\;.
\]
In order to establish this inequality we distinguish between two cases:
 $|x_1-y_1| \geq |y_1|/2$ and  $|x_1-y_1| \leq |y_1|/2$.

If $|x_1-y_1| \ge |y_1|/2$    then
$d_\delta(x_1,0\,;y_1,0)\leq 2\, d_\delta(0,0\,;y_1,0)$.
But $d_\delta(0,0\,;y_1,0)$ is less than the length of a straight line path
from $0$ to $y_1$.
Thus 
\[
d_\delta(x_1,0\,;y_1,0)\leq a\,
\frac{|y_1|^{\phantom{(\delta_1,\delta_1')}}}{|y_1|^{(\delta_1,\delta_1')}  }\leq
a\,\frac{|x_1-y_1|^{\phantom{(\delta_1,\delta_1')} }}{(|x_1|+|y_1|)^{(\delta_1,\delta_1')}  }
\;\;\;.
\]
If, however,   $|x_1-y_1| \leq |y_1|/2$ then    we consider the path
$(x_1(t),x_2(t))=(tx_1+(1-t)y_1,0)$. 
Note that
 $|tx_1+(1-t)y_1| \le |y_1|-t |x_1-y_1|  \leq  |y_2|/2$
so
\[
d_\delta(x_1,0\,;y_1,0)\leq a \int_0^1dt\, |x_1-y_1| \,c_{\delta_1,\delta_1'}(|y_1|/2)^{-1/2}
\sim \frac{|x_1-y_1|^{\phantom{(\delta_1,\delta_1')} }}{(|x_1|+|y_1|)^{(\delta_1,\delta_1')}  }
\;\;\;.
\]
This completes the bound of the first term on the right hand side of (\ref{egr5.1}).

Next we bound the second term on the right  of (\ref{egr5.1}).
Specifically we will establish that
\begin{equation}
 d_\delta(y_1,0\,;y_1,y_2-x_2) \le a\, \Delta_\delta (y_1,0\,;y_1,x_2-y_2)\sim
\Delta_\delta (x_1,x_2\,;y_1,y_2)
\;\;\;.\label{csg.one}
\end{equation}
If $|x_2-y_2|\le |y_1|^{(\rho,\rho')} $
then considering the path $y(t)=(y_1,t(x_2-y_2))$ we find
\begin{eqnarray}
 d_\delta(y_1,0\,;y_1,y_2-x_2)&\leq&
a\int_0^1dt\,  |x_2-y_2|\,c_{\delta_2,\delta_2'}(|y_1|)^{-1/2}\nonumber\\[5pt]
&\sim &{|x_2-y_2|}{|y_1|^{(-\delta_2,-\delta_2')} }
\sim \Delta_\delta (x_1,x_2\,;y_1,y_2)
\;\;\;.
\label{ecsg.a1}
\end{eqnarray}
If, however,
$|x_2-y_2|\geq |y_1|^{(\rho,\rho')} $ we set $\tilde{y}_1=
(y_1/{|y_1|})|x_2-y_2|^{(1/\rho,1/\rho')}$.
Then by (\ref{ecsg.a1})
\[
d_\delta(\tilde{y}_1,0\,;\tilde{y}_1, y_2-x_2)\leq
a\,{|x_2-y_2|}|\tilde{y}_1|^{(-\delta_2,-\delta_2')}  \sim
{|x_2-y_2|^{(1-\gamma,1-\gamma')}}
\;\;\;.
\]
Therefore
\begin{eqnarray*}
 d_\delta(y_1,0\,;y_1,y_2-x_2)&\leq&
2 \,d_\delta(\tilde{y}_1,0;\tilde{y}_1,0)+
d_\delta(\tilde{y}_1,0;\tilde{y}_1, y_2-x_2)\\[5pt]
&\leq& a\, ( |\tilde{y}_1|^{(1-\delta_1,1-\delta_1')}+
{|x_2-y_2|^{(1-\gamma,1-\gamma')}})\\[5pt]
&\sim& {|x_2-y_2|^{(1-\gamma,1-\gamma')}}
\sim  \Delta_\delta (x_1,x_2\,;y_1,y_2)
\;\;\;.
\end{eqnarray*}
Combination of these estimates then gives an upper bound
\[
d_\delta(x_1,x_2\,;y_1,y_2)\leq a\, D_\delta (x_1,x_2\,;y_1,y_2)
\]
for all $x,y$.
Therefore to complete the proof of equivalence of the distances
we have to establish a similar lower bound.

\bigskip

First we argue that
\begin{equation}
  c_{\delta_1,\delta_1'}(x_1)(\nabla_{x_1}D_\delta (x_1,x_2;y_1,y_2))^2
+c_{\delta_2,\delta_2'}(x_1)(\nabla_{x_2}D_\delta (x_1,x_2;y_1,y_2))^2
\leq a
\;\;\;.\label{csg.nabla}
\end{equation}
To establish (\ref{csg.nabla}) we first note that
$\nabla_{x_2}({|x_1-y_1|}/{(|x_1|+|y_1|)^{(\delta_1,\delta_1')}})=0$.
Secondly,
\begin{eqnarray*}
\nabla_{x_1}\frac{|x_1-y_1|^{\phantom{(\delta_1,\delta_1')}}}
{(|x_1|+|y_1|)^{(\delta_1,\delta_1')}}=
({\nabla_{x_1}|x_1-y_1|})\,{(|x_1|+|y_1|)^{ (-\delta_1,-\delta_1') }}+
|x_1-y_1|\,\nabla_{x_1}(|x_1|+|y_1|)^{(-\delta_1,-\delta_1')}\\[5pt]
 \le  (|x_1|+|y_1|)^{(-\delta_1,-\delta_1') }  +a\,
\frac{|x_1-y_1|^{\phantom{(1+\delta_1,1+\delta_1')}}}
{(|x_1|+|y_1|)^{(1+\delta_1,1+\delta_1')}}
\le a\,(|x_1|+|y_1|)^{(-\delta_1,-\delta_1')}\;\;\;.
\end{eqnarray*}
Therefore
\[
 c_{\delta_1,\delta_1'}(x_1)\left(\nabla_{x_1}
\frac{|x_1-y_1|^{\phantom{(\delta_1,\delta_1')}}}{(|x_1|+|y_1|)^{(\delta_1,\delta_1')}}\right)^2
+c_{\delta_2,\delta_2'}(x_1)\left(\nabla_{x_2}
\frac{|x_1-y_1|^{\phantom{(\delta_1,\delta_1')}}}{(|x_1|+|y_1|)^{(\delta_1,\delta_1')}}\right)^2 \le a
\;\;\;.
\]
Hence to establish  (\ref{csg.nabla}) it is enough to show that
\begin{equation}\label{csg.nabla1}
  c_{\delta_1,\delta_1'}(x_1)(\nabla_{x_1}\Delta_\delta
(x_1,x_2;y_1,y_2))^2
+c_{\delta_2,\delta_2'}(x_1)(\nabla_{x_2}\Delta_\delta
(x_1,x_2;y_1,y_2))^2
\le a
\;\;\;.
\end{equation}
Now
\[
c_{\delta_1,\delta_1'}(x_1)(\nabla_{x_1}|x_2-y_2|^{(1-\gamma,1-\gamma')})^2=0
\]
and
\[
c_{\delta_2,\delta_2'}(x_1)(\nabla_{x_2}|x_2-y_2|^{(1-\gamma,1-\gamma')})^2
\leq \left(\frac{|x_1|^{(\delta_2,\delta_2')}}{
|x_2-y_2|^{(\gamma,\gamma')}}\right)^2
\le a
\]
for all  $|x_2-y_2|\ge |x_1|+|y_1|^{(\rho,\rho')} $.
Next
\begin{eqnarray*}
   c_{\delta_1,\delta_1'}(x_1) \left(\nabla_{x_1}
\frac{|x_2-y_2|^{\phantom{(\delta_2,\delta_2')}} }{  (|x_1|+|y_1|)^{(\delta_2,\delta_2')}   }\right)^2
&\leq& a\,  c_{\delta_1,\delta_1'}(x_1)
\left(\frac{|x_2-y_2|^{ \phantom{(1+\delta_2,1+\delta_2')} }}{(|x_1|+|y_1|)^{ (1+\delta_2,1+\delta_2') }
}\right)^2\\[5pt]
&&\hspace{-2cm}{}\leq a\,
\left(\frac{|x_2-y_2|\,|x_1|^{ (\delta_1,\delta_1')\phantom{1} } }{ (|x_1|+|y_1|)^{
(1+\delta_2,1+\delta_2') }}\right)^2 \leq a\,\left(\frac{|x_2-y_2|^{\phantom{ (\rho,\rho')} }}{
(|x_1|+|y_1|)^{ (\rho,\rho') }}\right)^2\le a
\end{eqnarray*}
for all  $|x_2-y_2|\le |x_1|+|y_1|^{(\rho,\rho')} $.
Finally
\[
 c_{\delta_2,\delta_2'}(x_2)\left(\nabla_{x_2}
\frac{|x_2-y_2|^{\phantom{(\delta_2,\delta_2')}} }{  (|x_1|+|y_1|)^{(\delta_2,\delta_2')}   }\right)^2
\leq \left(
\frac{|x_1|^{(\delta_2,\delta_2')}  }{
(|x_1|+|y_1|)^{(\delta_2,\delta_2')}   }\right)^2
\le a
\;\;\;.
\]
This completes the verification of (\ref{csg.nabla}).

  Since $D_\delta$ satisfies (\ref{csg.nabla}) it follows formally from the definition (\ref{ecsg1.21})
  of the Riemannian distance that $D_\delta (x_1,x_2\,;y_1,y_2) \leq a\, d_\delta (x_1,x_2\,;y_1,y_2)$.
 One cannot, however, immediately make this deduction since $D_\delta\not\in W^{1,\infty}(\Ri^n\times \Ri^m)$. 
  It is, however, a  continuous    function which is locally Lipschitz differentiable on
$(\Ri^n -   \{ 0  \} )\times \Ri^m$.
But the distance is not changed if one replaces the space of trial functions
$ W^{1,\infty}(\Ri^n\times \Ri^m)$ in the definition  (\ref{ecsg1.21}) by a space of functions which are Lipschitz differentiable on the complement  of a closed set of measure zero.
This can be deduced by remarking that since we may assume the coefficients are continuous both 
definitions agree with the shortest path definition of the distance. 
Therefore $D_\delta (x_1,x_2\,;y_1,y_2) \leq a\, d_\delta (x_1,x_2\,;y_1,y_2)$ and since we have already established the converse inequality one concludes that
$D_\delta (x_1,x_2\,;y_1,y_2) \sim d_\delta (x_1,x_2\,;y_1,y_2)$.

\smallskip

It remains to  prove the volume estimates (\ref{evol}).
The proof will be divided into three steps.
First we consider small $r$, small compared with $|x_1|$, secondly we consider 
large $r$ and finally we deal with  intermediate values.

\smallskip

\noindent{\bf Step 1}$\;$ Assume $c\,r \leq |x_1|^{(1-\delta_1,1-\delta_1')}$ with $c>1$.
In fact we will choose $c\gg1$ in the course of the proof.
We now argue that there are $a_1, a_2\in\Ri$ with $0<a_1<a_2$ such that 
\begin{eqnarray}
(x_1,x_2)&&+
\left[-a_1\,r\,|x_1|^{(\delta_1,\delta_1')}, a_1\,r\,|x_1|^{(\delta_1,\delta_1')}\right]^n
\times
\left[-a_1\,r\,|x_1|^{(\delta_2,\delta_2')}, a_1\,r\,|x_1|^{(\delta_2,\delta_2')}\right]^m
\subset  B(x_1,x_2\,;r)\nonumber \\[5pt]
&&\hspace{-1cm}{}\subset
(x_1,x_2)+
\left[-a_2\,r\,|x_1|^{(\delta_1,\delta_1')}, a_2\,r\,|x_1|^{(\delta_1,\delta_1')}\right]^n
\times
\left[-a_2\,r\,|x_1|^{(\delta_2,\delta_2')}, a_2\,r\,|x_1|^{(\delta_2,\delta_2')}\right]^m\;\;.
\label{incl1}
\end{eqnarray}
Once this is established one has the volume estimates
\begin{equation}
|B(x_1,x_2\,;r)|\sim r^{n+m}|x_1|^{(n\delta_1+m\delta_2,n\delta_1'+m\delta_2')}
=r^{n+m}|x_1|^{(\beta,\beta')}
\label{ecsg5.101}
\end{equation}
for $r \leq c^{-1}\, |x_1|^{(1-\delta_1,1-\delta_1')}$

First consider the right hand inclusion of (\ref{incl1}).
Set $\|x\|=\max_{1\leq k\leq n}|x^{(k)}|$ and $\|y\|=\max_{1\leq l\leq m}|y^{(l)}|$
where $x^{(k)}$ and $y^{(l)}$ are the components of $x\in\Ri^n$ and $y\in \Ri^m$, respectively.
Thus we  have to prove that if $(y_1,y_2)\in B(x_1,x_2\,;r)$ then 
$\|x_1-y_1\|\leq a_2\,r\,|x_1|^{(\delta_1,\delta_1')}$ and 
$\|x_2-y_2\|\leq a_2\,r\,|x_1|^{(\delta_2,\delta_2')}$.
But $d_\delta (x_1,x_2;y_1,y_2)\sim D_\delta (x_1,x_2;y_1,y_2)$.
Hence if $(y_1,y_2)\in B(x_1,x_2\,;r)$ then 
$\|x_1-y_1\|\leq |x_1-y_1|\leq a\,r\,(|x_1|+|y_1|)^{(\delta_1,\delta_1')}$.
Therefore
\begin{eqnarray*}
|\|x_1\|-\|y_1\||\leq \|x_1-y_1\|&\leq &a'\,c^{-1}\,(|x_1|+|y_1|)^{(1-\delta_1,1-\delta_1')}\,(|x_1|+|y_1|)^{(\delta_1,\delta_1')}\\[5pt]
&=&a'\,c^{-1}\,(|x_1|+|y_1|)\leq a'\,c^{-1}\,(n+m)(\|x_1\|+\|y_1\|)
\;\;\;.
\end{eqnarray*}
Choosing $c$ large one deduces that $\|y_1\|\sim\|x_1\|$ and $|y_1|\sim |x_1|$.
In particular one has 
\[
\|x_1-y_1\|\leq a\,r\,(|x_1|+|y_1|)^{(\delta_1,\delta_1')}\leq a_2\,r\,|x_1|^{(\delta_1,\delta_1')}
\]
as required.

Secondly, since  $(y_1,y_2)\in B(x_1,x_2\,;r)$ one has $ \Delta_\delta (x_1,x_2\,;y_1,y_2)\leq a\,r$.
There are two cases to consider.  
The first is if $ |x_2-y_2|\geq (|x_1|+|y_1|) ^{(\rho,\rho')}$  then   $\Delta_\delta (x_1,x_2\,;y_1,y_2)= 
|x_2-y_2|^{(1-\gamma,1-\gamma')}  $.
Hence 
\[
|x_2-y_2|^{(1-\gamma,1-\gamma')}  \leq a\,r\leq a\,c^{-1}\,|x_1|^{(1-\delta_1,1-\delta_1')}
\leq a\,c^{-1}\,(|x_1|+|y_1|)^{(1-\delta_1,1-\delta_1')}
\;\;\;.
\]
Thus if $c\geq a$ then 
\[
|x_2-y_2|\leq\,(|x_1|+|y_1|)^{((1-\delta_1)(1-\gamma)^{-1},(1-\delta_1')(1-\gamma')^{-1})}
=(|x_1|+|y_1|)^{(\rho,\rho')}
\]
which is in contradiction with the assumption  $ |x_2-y_2|\geq (|x_1|+|y_1|) ^{(\rho,\rho')}$.
Therefore one must have  $ |x_2-y_2|\leq (|x_1|+|y_1|) ^{(\rho,\rho')}$ and 
$\Delta_\delta (x_1,x_2\,;y_1,y_2)= {|x_2-y_2|}/{(|x_1|+|y_1|)^{(\delta_2,\delta_2')}}$.
Then, however,
\[
\|x_2-y_2\|\leq |x_2-y_2|\leq a\,r\,(|x_1|+|y_1|)^{(\delta_2,\delta_2')}\leq a_2\,r\,|x_1|^{(\delta_2,\delta_2')}
\]
because $|y_1|\sim|x_1|$ by the previous argument.
Thus the proof of the right hand inclusion of (\ref{incl1}) is complete.

Secondly, consider the left hand inclusion of (\ref{incl1}).
Now we need to prove that if one has 
$\|x_1-y_1\|\leq a_1\,r\,|x_1|^{(\delta_1,\delta_1')}$ and 
$\|x_2-y_2\|\leq a_1\,r\,|x_1|^{(\delta_2,\delta_2')}$ then 
$(y_1,y_2)\in B(x_1,x_2\,;r)$.
But by the first assumption
\[
|x_1-y_1|/(|x_1|+|y_1|)^{(\delta_1,\delta_1')}\leq  a_1\,n\,r\,|x_1|^{(\delta_1,\delta_1')}/(|x_1|+|y_1|)^{(\delta_1,\delta_1')}\leq a_1\,n\,r
\;\;\;.
\]
Then by the second assumption
\begin{eqnarray*}
|x_2-y_2|\leq a_1\,m\,r\,|x_1|^{(\delta_2,\delta_2')}
&\leq& a_1\,m\,c^{-1}\,|x_1|^{(1+\delta_2-\delta_1,1+\delta_2'-\delta_1')}\\[5pt]
&=&a_1\,m\,c^{-1}\,|x_1|^{(\rho,\rho')}
\leq a_1\,m\,c^{-1}\,(|x_1|+|y_1|)^{(\rho,\rho')}
\;\;\;.
\end{eqnarray*}
Hence if $c\geq a_1\,m$ then $|x_2-y_2|\leq  (|x_1|+|y_1|)^{(\rho,\rho')}$.
Therefore it follows that  $\Delta_\delta (x_1,x_2\,;y_1,y_2)= {|x_2-y_2|}/{(|x_1|+|y_1|)^{(\delta_2,\delta_2')}}$.
Hence using $\|x_2-y_2\|\leq a_1\,r\,|x_1|^{(\delta_2,\delta_2')}$  again one has
\[
\Delta_\delta (x_1,x_2\,;y_1,y_2)\leq  a_1\,m\,r\,|x_1|^{(\delta_2,\delta_2')}/{(|x_1|+|y_1|)^{(\delta_2,\delta_2')}}\leq a_1\,m \,r
\]
and 
\[
D_\delta (x_1,x_2\,;y_1,y_2)=
{|x_1-y_1|}/{(|x_1|+|y_1|)^{(\delta_1,\delta_1')}}
+\Delta_\delta (x_1,x_2\,;y_1,y_2)\leq a_1\,(n+m)\,r
\;\;\;.
\]
Hence if $a_1$ is sufficiently small one concludes that $(y_1,y_2)\in B(x_1,x_2\,;r)$.

\smallskip

\noindent{\bf Step 2}$\;$ Assume $r/c \geq |x_1|^{(1-\delta_1,1-\delta_1')}$
with $c>1$ where we will again choose $c\gg1$.
We now argue that there are
$a_1, a_2\in\Ri$ with $0<a_1<a_2$ such that 
\begin{eqnarray}
(0,x_2)+
\left[-a_1\,r^{(\sigma,\sigma')},
a_1\,r ^{(\sigma,\sigma')}\right]^n
&\times&
\left[-a_1\,r ^{(\rho\sigma,\rho'\sigma')},
a_1\,r ^{(\rho\sigma,\rho'\sigma')} \right]^m
\subset B(x_1,x_2\,;r)\nonumber\\[5pt]
&&\hspace{-4.5cm}{}\subset
(0,x_2)+
\left[-a_2\,r^{(\sigma,\sigma')},
a_2\,r ^{(\sigma,\sigma')}\right]^n
\times
\left[-a_2\,r ^{(\rho\sigma,\rho'\sigma')},
a_2\,r ^{(\rho\sigma,\rho'\sigma')}\right]^m
\label{incl2}
\end{eqnarray}
where $\sigma=(1-\delta_1)^{-1}$ and $\sigma'=(1-\delta'_1)^{-1}$.
These inclusions then yield the volume estimates
\begin{equation}
|B(x_1,x_2\,;r)|\sim r^{(n\sigma, n\sigma')}r^{(m\rho\sigma, m\rho'\sigma')}=r^{(D,D')}
\;\;\;.
\label{ecsg5.102}
\end{equation}

The first step in deducing the inclusions (\ref{incl2}) is to observe that if $c$ is sufficiently large
then
\[
B(0,x_2\,;r/2)\subseteq B(x_1,x_2\,;r)\subseteq B(x_1,x_2\,;2r)
\]
and so the proof is effectively reduced to the case $x_1=0$.
But then the rest of the proof is similar to the argument given in Step~1 but somewhat simpler
because of the choice $x_1=0$.
We omit the details.

\smallskip

\noindent{\bf Step 3}$\;$
It follows from Step~1, and in particular (\ref{ecsg5.101}), that the volume estimates (\ref{evol}) are valid
if $ |x_1|^{(1-\delta_1,1-\delta_1')}\geq c\, r$.
Alternatively, it follows from Step~2, and in particular (\ref{ecsg5.102}),
that the estimates  are valid if  $r \geq c\, |x_1|^{(1-\delta_1,1-\delta_1')}$.
Therefore we now assume that $r\in[c^{-1} |x_1|^{(1-\delta_1,1-\delta_1')},c \,|x_1|^{(1-\delta_1,1-\delta_1')}]$.
But then setting $r_1=r/c$ and $r_2=r\,c$ one has
 $B(x_1,x_2\,;r_1) \subset B(x_1,x_2\,;r) \subset B(x_1,x_2\,;r_2)$.
 But $r_1\leq  |x_1|^{(1-\delta_1,1-\delta_1')}$.
 Hence
 $B(x_1,x_2\,;r_1)~\sim r_1^{n+m}|x_1|^{(\beta,\beta')}\sim
 r^{n+m}\,|x_1|^{(\beta,\beta')}$ by Step~1.
 Moreover, $r_2\geq|x_1|^{(1-\delta_1,1-\delta_1')}$.
 Hence 
$B(x_1,x_2\,;r_2)~\sim r_2^{(D,D')} \sim
 r^{(D,D')}$
 by Step~2.
 Combining these estimates one concludes that 
 \[
 a\,r^{n+m}\,|x_1|^{(\beta,\beta')}\leq |B(x_1,x_2\,;r) |\leq a'\,r^{(D,D')}
 \]
 for all $r\in[c^{-1} |x_1|^{(1-\delta_1,1-\delta_1')},c \,|x_1|^{(1-\delta_1,1-\delta_1')}]$.
But in this range $r^{(D,D')}\sim a\,r^{n+m}\,|x_1|^{(\beta,\beta')}$.
Therefore  the volume estimates are established again and the proof of the proposition is complete.
 \hfill$\Box$
 
 \ruimte
 
 The volume estimates allow one to prove the doubling property and to identify the doubling dimension.
 
 \begin{cor}\label{ccsg5.1} The Riemannian balls $B(x_1,x_2\,;r)$ associated with the Gru\v{s}in operator satisfy the doubling property
 \[
 |B(x_1,x_2\,;s\,r)|\leq a\,s^{(D\vee D')} |B(x_1,x_2\,;r)|
 \]
 for all $(x_1,x_2)\in \Ri^{n+m}$ and all $s\geq 1$.
 \end{cor}
 \proof\
 There are three cases to consider.
 
 \smallskip
 
\noindent{\bf Case~1}$\;\;r\leq r\,s\leq |x_1|^{(1-\delta_1,1-\delta_1')}$.
Then the volume estimates of Proposition~\ref{pcsg5.1} give
 \[
 |B(x_1,x_2\,;s\,r)|\sim(r\,s)^{n+m}\,|x_1|^{(\beta,\beta')}\sim s^{n+m}\, |B(x_1,x_2\,;r)|
 \;\;\;.
 \]
 But $n+m\leq D\vee D'$ and $s\geq 1$ so $s^{n+m}\leq s^{(D\vee D')}$ and the doubling property follows.

\smallskip

\noindent{\bf Case~2}$\;\;|x_1|^{(1-\delta_1,1-\delta_1')}\leq r\leq rs$.
Then the volume estimates  give
 \[
 |B(x_1,x_2\,;s\,r)|\sim(r\,s)^{(D,D')}\leq s^{(D\vee D')}\, r^{(D,D')}\sim s^{(D\vee D')}\,  |B(x_1,x_2\,;r)|
\]
  and the doubling property is established.
  
\smallskip

\noindent{\bf Case~3}$\;\;r\leq |x_1|^{(1-\delta_1,1-\delta_1')}\leq  rs$.
Then the volume estimates  give
 \[
 |B(x_1,x_2\,;s\,r)|\sim(r\,s)^{(D,D')}\leq s^{(D\vee D')}\, r^{(D,D')}
 \;\;\;.
\]
But
\[
r^{(D,D')}=r^{n+m}r^{(\beta(1-\delta_1)^{-1},\beta'(1-\delta_1')^{-1})}\leq 
r^{n+m}|x_1|^{(\beta,\beta')}\sim |B(x_1,x_2\,;r)|
\;\;\;.
\]
Combination of these estimates gives the doubling property again.
 \hfill$\Box$

\section{Kernel bounds}\label{Scsg6} 

In this section we derive two basic properties of the semigroup kernel $K_t$ associated with a
general Gru\v{s}in  operator $H$. 
Initially we show that the kernel conserves probability or, in probabilistic terms, that it is stocastichally complete.
Then  we establish that it satisfies off-diagonal volume dependent bounds.
Subsequently it is  possible to apply standard reasoning to obtain further more detailed properties 
such as Gaussian bounds and on diagonal lower bounds.
This will be discussed at the end of the section.

First, since the closure of the form $h$ associated with the Gru\v{s}in operator $H$ is a Dirichlet form the corresponding semigroup $S$ is submarkovian. 
In particular it extends to a contractive 
semigroup on $L_\infty(\Ri^n\times\Ri^m)$.
In terms of the kernel this means that 
\[
0\leq \esssup_{x\in\Ri^{n+m}}\int_{\Ri^{n+m}}dy\,K_t(x\,;y)\leq 1
\;\;\;.
\]
But the particular structure of the operator gives a stronger result.

\begin{thm}\label{tcsg6.1} 
The semigroup $S$ associated with the Gru\v{s}in operator  on $L_\infty(\Ri^n\times\Ri^m)$ satisfies
$S_t\one=\one$ for all $t>0$.
Hence the kernel satisfies
\begin{equation}
\int_{\Ri^n\times\Ri^m}dy\,K_t(x\,;y)=1
\label{ecsg6.1}
\end{equation}
for  all $t>0$ almost all $x\in \Ri^{n+m}$.
\end{thm}
\proof\
Given that the semigroup satisfies the $L_2$ off-diagonal bounds of Proposition~\ref{pscsg4.0}
and the volume of the Riemannian balls have polynomial growth, by the estimates of 
Proposition~\ref{pcsg5.1}, one can prove the theorem by  a slight variation of the argument given in Proposition~3.6 of \cite{ERSZ1} but estimating with respect to a Riemannian distance instead of the Euclidean distance.
We will, however, give a different variation on the argument which establishes a useful convergence
result.

Let  $H_{N,\varepsilon}$ be  the strongly
elliptic approximants,  with coefficients $C_{N,\varepsilon}=(C\wedge NI)+\varepsilon I$, to the Gru\v{s}in operator.
Then the  semigroups $S^{(N,\varepsilon)}_t$ generated by the   $H_{N,\varepsilon}$
converge strongly on $L_2$ to the semigroup $S_t$ in the double limit $N\to\infty$ followed by $\varepsilon\to 0$.
But the convergence is stronger.

\begin{prop}\label{pcsg4.9}
The semigroups $S^{(N,\varepsilon)}_t$ 
converge strongly to  $S_t$ on each of the  $L_p$-spaces with $p\in[1,\infty\rangle$
and in the weak$^*$ sense on $L_\infty$.
\end{prop}
\proof\
It suffices to prove the convergence on $L_1$.
Then weak$^*$ convergence on $L_\infty$ follows by duality and $L_p$-convergence 
for $p\in\langle1,\infty\rangle$ follows since 
\[
 \|(S^{(N,\varepsilon)}_t-S_t)\psi\|_p
 \leq  \|(S^{(N,\varepsilon)}_t-S_t)\psi\|_1^{1/p}\,
 (2\,\|\psi\|_\infty)^{1-1/p}
 \]
 for all $\psi$ in the dense subset $ L_1\cap L_\infty$.

 Since the semigroups are contractive it suffices to establish the $L_1$-convergence
 on a  subset of $L_1$ whose span is dense.
 In particular it suffices to prove convergence on $L_1(A)\cap L_2(A)$ for each bounded open
 subset $A$.
 Moreover one can restrict to positive functions.

 Fix $\varphi_A\in L_1(A)\cap L_2(A)$ and assume $\varphi_A$ is positive.
 Next let $B\supset A$ be  a bounded closed set.
Then
\begin{eqnarray*}
\|(S^{(N,\varepsilon)}_t-S_t)\varphi_A\|_1&\leq &
\|\one_B(S^{(N,\varepsilon)}_t-S_t)\varphi_A\|_1
+\|\one_{B^{\rm c}}S^{(N,\varepsilon)}_t\varphi_A\|_1+\|\one_{B^{\rm c}}S_t\varphi_A\|_1\\[5pt]
&\leq&|B|^{1/2}\|(S^{(N,\varepsilon)}_t-S_t)\varphi_A\|_2
+(\one_{B^{\rm c}},S^{(N,\varepsilon)}_t\varphi_A)+(\one_{B^{\rm c}},S_t\varphi_A)
\end{eqnarray*}
where we have used the positivity of the semigroups and the functions to express the norms as pairings between $L_1$ and $L_\infty$.

Next let 
 $d_{N,\varepsilon}$ denote the Riemannian distance corresponding to the coefficients 
 $C_{N,\varepsilon}$ and $d_1$ the distance corresponding to the coefficients $C+I$ then 
 $d_1\leq d_{N,\varepsilon}$ for all $N\geq 0$ and all $\varepsilon\in\langle0,1]$.
 But it  follows from Proposition~\ref{pscsg4.0} applied to $S^{(N,\varepsilon)}_t$ that 
 \[
|(\varphi_C,S^{(N,\varepsilon)}_t\varphi_A)|\leq e^{-d_{N,\varepsilon}(A ;C)^2 (4t)^{-1}}\|\varphi_A\|_2\|\varphi_C\|_2\leq e^{-d_1(A ;C)^2 (4t)^{-1}}\|\varphi_A\|_2\|\varphi_C\|_2
\]
for all open sets $C$ and all $\varphi_C\in L_2(C)$.
Now choose $R$ sufficiently large that $A\subseteq B_{1,R}=\{x:d_1(0\,;x)<R\}$ and let 
$B=\overline{B_{1,2R}}$. 
Then one can separate $B^{\rm c}$ into annuli and make a quadrature estimate, as in the proof of Proposition~3.6 of \cite{ERSZ1}, to find
\begin{eqnarray*}
(\one_{B^{\rm c}},S^{(N,\varepsilon)}_t\varphi_A)& \leq & \sum_{n\geq2}e^{-d_1(B_{1,(n+1)R} \backslash B_{1,nR} ;B_{1,R})^2 (4t)^{-1}}
     |B_{1, (n+1)R}|^{1/2} \|\varphi_A\|_2 
  \end{eqnarray*}
  uniformly for $N\geq 1$ and $\varepsilon\in\langle0,1]$.

But $d_1(B_{1,(n+1)R} \backslash B_{1,nR}\, ;B_{1,R})\geq (n-1)\,R$ by the triangle inequality.
Moreover $H_1$, the operator with coefficents $C+I$,
 is a Gru\v{s}in  operator with $\delta_1=0=\delta_2$.
 Hence   it follows that 
$ |B_{1, (n+1)R}|^{1/2}\leq a\, (nR)^{D'/2}$ for all $n\geq 2$ and $R\geq 1$
by Proposition~\ref{pcsg5.1}.
Therefore one obtains an estimate
\[
(\one_{B^{\rm c}},S^{(N,\varepsilon)}_t\varphi_A)\leq a\,\sum_{n\geq 2}(nR)^{D'/2}e^{-a'n^2R^2t^{-1}}
\|\varphi_A\|_2
\]
for all  $R\geq 1$, uniform in $N$ and $\varepsilon$.
A similar bound is valid for $S_t$ by the same reasoning.

Combining these estimates gives
\begin{eqnarray*}
\|(S^{(N,\varepsilon)}_t-S_t)\varphi_A\|_1&\leq &
|B_{1,2R}|^{1/2}\,\|(S^{(N,\varepsilon)}_t-S_t) \varphi_A\|_2
+2\,a_R\|\varphi_A\|_2
\end{eqnarray*}
where  $a_R\to 0$ as $R\to\infty$ uniformly in $N$ and $\varepsilon$.
Therefore 
\[
\limsup_{\varepsilon\to0,N\to\infty}\|(S^{(N,\varepsilon)}_t-S_t)\varphi_A\|_1
\leq 2\,a_R\|\varphi_A\|_2
\]
by the $L_2$ convergence.
Then the $L_1$ convergence follows because $a_R$ can be made arbitrarily small
by choosing $R$ sufficiently large.
 \hfill$\Box$

\bigskip

\noindent{\bf Proof of Theorem~\ref{tcsg6.1} }
Since the approximants $H_{N,\varepsilon}$ are strongly elliptic $S^{(N,\varepsilon)}_t\one=\one$
for all $t>0$.
Then by weak$^*$ convergence $S_t\one=\one$.
\hfill$\Box$

\bigskip

\begin{remarkn}\label{rcsg6.10}
The proof of Proposition~\ref{pcsg4.9} and Theorem~\ref{tcsg6.1} 
 uses very little structure of the Gru\v{s}in operator.
 It only requires $L_2$ off-diagonal bounds and polynomial volume growth.
 The first are given for general elliptic operators by Proposition~\ref{pscsg4.0} and the polynomial growth
 follows from the Gru\v{s}in structure.
 \end{remarkn}

Next we consider upper bounds on the semigroup kernel.

\begin{thm}\label{tcsg6.2} 
There is an $a>0$ such that the  semigroup kernel $K$ of the Gru\v{s}in operator $ H$
satisfies 
 \begin{equation}
0\leq  K_t(x\,;y)\leq a\,(|B(x\,;t^{1/2})|\,|B(y\,;t^{1/2})|)^{-1/2}
\label{ecsg6.2}
\end{equation}
for  all $t>0$ and  almost all $x,y\in \Ri^{n+m}$.
\end{thm}

Since the semigroup $S$ is self-adjoint it follows that the semigroup kernel is positive-definite.
Therefore one formally has $|K_t(x\,;y)|^2\leq K_t(x\,;x)K_t(y\,;y)$ and the estimate
apparently reduces to an on-diagonal estimate.
But this is only a formal calculation since the kernel is not necessarily continuous and its diagonal value
is not necessarily defined.
Nevertheless the starting point of the proof is a set-theoretic reduction to an on-diagonal estimate.

\begin{lemma}\label{lcsg6.1}
Let $X,Y$ be open sets and define $K_t(X\,;Y)=\esssup_{x\in X,y\in Y}K_t(x\,;y)$.
Then
\[
|K_t(X\,;Y)|^2\leq K_t(X\,;X)\,K_t(Y\,;Y)
\;\;\;.
\]
\end{lemma}
\proof\
First observe that 
\begin{eqnarray*}
K_t(X\,;Y)=\|\one_YS_t\one_X\|_{1\to\infty}
\leq\|S_{t/2}\one_X\|_{1\to2}\|\one_Y S_{t/2}\|_{2\to\infty}
= \|\one_XS_{t/2}\|_{2\to\infty}\|\one_Y S_{t/2}\|_{2\to\infty}
\;\;\;.
\end{eqnarray*}
But if $T$ is  bounded from $L_2$ to $L_\infty$ then $(\|T\|_{2\to\infty})^2=\|TT^*\|_{1\to\infty}$.
Therefore
\begin{eqnarray*}
|K_t(X\,;Y)|^2
\leq\|\one_XS_{t}\one_X\|_{1\to\infty}\|\one_Y S_{t}\one_Y\|_{1\to\infty}
=K_t(X\,;X)\,K_t(Y\,;Y)
\end{eqnarray*}
as required.\hfill$\Box$

\bigskip

\noindent{\bf Proof of Theorem~\ref{tcsg6.2} } It follows from the lemma that 
\[
K_t(x\,;y)\leq \inf_{X\ni x}K_t(X\,;X)^{1/2}\,\inf_{Y\ni y}K_t(Y\,;Y)^{1/2}
\;\;\;.
\]
Thus it suffices to prove that 
\[
 \inf_{X\ni x}K_t(X\,;X)\leq a\,|B(x\,;t^{1/2})|^{-1}
 \;\;\;.
 \]
 There are two distinct cases corresponding to the different volume behaviours given by 
 Proposition~\ref{pcsg5.1}.
 
 First, let  $x=(x_1,x_2)$ and suppose $|x_1|^{(1-\delta_1,1-\delta_1')}\leq t^{1/2}$. 
 Then 
 $|B(x\,;t^{1/2})|\sim t^{(D/2,D'/2)}$ by  Proposition~\ref{pcsg5.1}.
 On the other hand $\|S_t\|_{1\to\infty}\leq a\,t^{(-D/2,-D'/2)}$ by Proposition~\ref{pcsg3.0}.
 Therefore 
 \[
 \inf_{X\ni x}K_t(X\,;X)\leq \|K_t\|_\infty=\|S_t\|_{1\to\infty}\leq a\,t^{(-D/2,-D'/2)}
 \leq a'\,|B(x\,;t^{1/2})|^{-1}
 \;\;\;.
 \]
 
Secondly, suppose that $t^{1/2}\leq |x_1|^{(1-\delta_1, 1-\delta_1')}$. 
 Then  $|B(x\,;t^{1/2})|\sim  t^{(n+m)/2}|x_1|^{(\beta,\beta')}$ by  Proposition~\ref{pcsg5.1}.
 This case is considerably more difficult to analyze and it is here that we apply the comparison techniques of Section~\ref{Scsg4}.
 
Set $r=|x_1|$.
Let $C$ denote the coefficient matrix of $H$ and  choose $a_n, a_m>0$ such that 
$C(y)\geq a_n\,|y_1|^{(2\delta_1,2\delta_1')} I_n+a_m\,|y_1|^{(2\delta_2,2\delta_2')} I_m$ for all $y=(y_1,y_2)$ with $|y_1|\leq r/2$.
Next set $C_r(y)=C(y)$ if $|y_1|>r/2$ and
$C_r(y)= a_n\,r^{(2\delta_1,2\delta_1')} I_n+a_m\,r^{(2\delta_2,2\delta_2')} I_m$ if $|y_1|\leq r/2$.
Then $C_r\geq a\,C$ for a suitable $a>0$. 
Let  $H_r$ be  the Gru\v{s}in operator with coefficient matrix $C_r$ 
and set $H_1=H_r$ and $H_2=a\,H$.
Then $H_1\geq H_2$.
Moreover $H_1\geq \mu I>0$  for some 
$\mu=a\,(a_nr^{(2\delta_1,2\delta_1')}\wedge a_m r^{(2\delta_2,2\delta_2')})$.
Thus the basic assumptions of Theorem~\ref{tcsg3.1} and Corollary~\ref{ccsg3.1} are satisfied with
this choice.
Now $H_r$ is a  Gru\v{s}in operator  with local parameters $\delta_1=0=\delta_2$ but with the same
global parameters $\delta_1',\delta_2'$ as $H$.
Therefore
\[
\|S^{(1,0)}_t\|_{1\to\infty}\vee \|S^{(2,0)}_t\|_{1\to\infty}\leq V(t)^{-1}
\]
for all $t>0$ where $V(t)=a\,t^{(D/2,D'/2)}$.
Since $V$ satisfies the doubling property one may now apply Corollary~\ref{ccsg3.1}.
Note that for this application $U=\{y: |y_1|\leq r/2\}$ and $d(x\,;U)\sim r^{(1-\delta_1,1-\delta_1')}$.

Next we need an improved estimate on the crossnorm $\|S^{(1,0)}_t\|_{1\to\infty}$ of the comparison 
semigroup.
Let $\widehat H$ denote the constant coefficient operator with coefficients $\widehat C_r
= a_n\,r^{(2\delta_1,2\delta_1')} I_n+a_m\,r^{(2\delta_2,2\delta_2')} I_m$.
Then $\widehat H$ is a Fourier multiplier and the corresponding function $F$  is given by
 $F(p_1,p_2)=a_n\,r^{(2\delta_1,2\delta_1')}p_1^2+a_m\,r^{(2\delta_2,2\delta_2')} p_2^2$.
But there is an $a>0$ such that $H_1 \geq a\, \widehat H$ and we can apply Lemma~\ref{lcsg2.2}
to obtain a uniform bound on $ \|S^{(1,0)}_t\|_{1\to\infty}$.
First, suppose $r=1$ then one immediately has $ \|S^{(1,0)}_t\|_{1\to\infty}\leq a\,t^{-(n+m)/2}$
for all $t>0$.
Secondly, the introduction of $r$ corresponds to a dilation of $\Ri^{n}\times \Ri^{m}$
with each direction in $\Ri^{n}$ dilated by $r^{(-\delta_1,-\delta_1')}$  and each direction in 
$\Ri^{m}$ dilated by $r^{(-\delta_2,-\delta_2')}$.
The dilation adds a factor to the crossnorm corresponding to the Jacobian
$r^{(-n\delta_1-m\delta_2,-n\delta_1'-m\delta_2')}=r^{(-\beta,-\beta')}$ of the dilation. 
Therefore  
\[
 \|S^{(1,0)}_t\|_{1\to\infty}\leq a'\,t^{-(n+m)/2}r^{(-\beta,-\beta')}
 \leq a''\,|B(x\,;t^{1/2})|^{-1}
 \]
 where the last bound follows from the second estimate of Proposition~\ref{pcsg5.1}.

Now we may apply Corollary~\ref{ccsg3.1} with $A=X$.
One obtains an estimate
 \begin{eqnarray*}
 \inf_{X\ni x}K_t(X\,;X)&=& \inf_{X\ni x}K^{(2,0)}_t(X\,;X)\\[5pt]
 &\leq&
 \|S^{(1,0)}_t\|_{1\to\infty}+
\mathop{\smash \inf\vphantom\sup}_{X\ni x}\sup_{y,z\in X}
 |K^{(1,0)}_t(y\,;z)-K^{(2,0)}_t(y\,;z)|\\[5pt]
 &\leq& a\,|B(x\,;t^{1/2})|^{-1}+a\,V(t^2/\rho^2)^{-1}\,(\rho^2/t)^{-1/2}
 \,e^{-\rho^2/(4t)}\\[5pt]
 &\leq& a\,|B(x\,;t^{1/2})|^{-1}\left(1+R(x\,;t)\right)
\end{eqnarray*}
where 
\[
R(x\,;t)=|B(x\,;t^{1/2})|\,V(t^2/\rho^2)^{-1}\,(\rho^2/t)^{-1/2}
 \,e^{-\rho^2/(4t)}
 \;\;\;,
 \]
 $\rho=d(x\,;U)\sim r^{(1-\delta_1,1-\delta_1')}$ and $V(t)=a\,t^{(D/2,D'/2)}$.
Thus it suffices to show that $R(x\,;t)$ is uniformly bounded for $x,t$ satisfying
 $t^{1/2}\leq r^{(1-\delta_1, 1-\delta_1')}$ with $r=|x_1|$. 
It is necessary to distinguish between three cases.

\smallskip

\noindent{\bf Case~1}$\;\;t,r\leq 1$.
In this case  $t^{1/2}\leq r^{1-\delta_1}\leq 1$ and $\rho\sim r^{1-\delta_1}$. 
Hence $\rho^2/t\sim (r^{1-\delta_1}t^{-1/2})^2$ and  $t^2/\rho^2\leq a\,t$. 
Therefore, since  $t\leq1$, one must have
$V(t^2/\rho^2)^{-1}\leq a\, t^{-D/2}$.
Moreover, $|B(x\,;t^{1/2})|\sim  t^{(n+m)/2}r^{\beta}$.
Then one  computes that 
\begin{eqnarray*}
R(x\,;t)&\leq&
  a\,t^{(n+m)/2}r^{\beta}\,t^{-D/2}\,(r^{1-\delta_1}t^{-1/2})^{-1}e^{-a'(r^{1-\delta_1}t^{-1/2})^2}\\[5pt]
  &=& a\,(r^{1-\delta_1}t^{-1/2})^{-1+\beta(1-\delta_1)^{-1}}\,e^{-a'(r^{1-\delta_1}t^{-1/2})^2}
  \;\;\;.
  \end{eqnarray*}
  But since $r^{1-\delta_1}t^{-1/2}\geq1$  it follows that $R(x\,;t)$ is uniformly bounded.

\smallskip

\noindent{\bf Case~2}$\;\;t,r\geq 1$.
In this case  $1\leq t^{1/2}\leq r^{1-\delta_1'}$ and $\rho\sim r^{1-\delta_1'}$. 
Hence $\rho^2/t\sim (r^{1-\delta_1'}t^{-1/2})^2$ and  $t^2/\rho^2\leq a\,t$. 
But now  $t\geq1$ and so 
$V(t^2/\rho^2)^{-1}\leq a\, t^{-D'/2}$.
Moreover, $|B(x\,;t^{1/2})|\sim  t^{(n+m)/2}r^{\beta'}$.
Then one  computes  as in Case~1 that
\begin{eqnarray*}
R(x\,;t)&\leq&
  a\,t^{(n+m)/2}r^{\beta'}\,t^{-D'/2}\,(r^{1-\delta_1'}t^{-1/2})^{-1}e^{-a'(r^{1-\delta_1'}t^{-1/2})^2}\\[5pt]
  &=& 
 a\,(r^{1-\delta_1'}t^{-1/2})^{-1+\beta'(1-\delta_1')^{-1}}\,e^{-a'(r^{1-\delta_1'}t^{-1/2})^2}
  \;\;\;.
  \end{eqnarray*}
  But since $r^{1-\delta_1'}t^{-1/2}\geq1$  it again  follows that $R(x\,;t)$ is uniformly bounded.

\smallskip

\noindent{\bf Case~3}$\;\;t\leq 1, r\geq 1$.
This is a hybrid case which is rather different to the previous two cases.
One again has $\rho\sim r^{1-\delta_1'}$ and  $\rho^2/t\sim (r^{1-\delta_1'}t^{-1/2})^2$.
Moreover, $t^2/\rho^2\leq a\,t$. 
But since $t\leq 1$ one has $V(t^2/\rho^2)^{-1}\leq a\, t^{-D/2}$.
In addition $|B(x\,;t^{1/2})|\sim  t^{(n+m)/2}r^{\beta'}$.
Therefore one now estimates that
\begin{eqnarray*}
R(x\,;t)&\leq&
  a\,t^{(n+m)/2}r^{\beta'}\,t^{-D/2}\,(r^{1-\delta_1'}t^{-1/2})^{-1}e^{-a'(r^{1-\delta_1'}t^{-1/2})^2}\\[5pt]
   &= &a\,r^{\beta'}\,(t^{-1/2})^{\beta(1-\delta_1)^{-1}}\,
  (r^{1-\delta_1'}t^{-1/2})^{-1} e^{-a'(r^{1-\delta_1'}t^{-1/2})^2}
  \;\;\;.
  \end{eqnarray*}
  But now $r^{1-\delta_1'}t^{-1/2}\geq 1$.
 Consequently for each $N\geq 1$ there is an $a_N$ such that 
\[
R(x\,;t)\leq 
a_N\,r^{\beta'}\,(t^{-1/2})^{\beta(1-\delta_1)^{-1}}\,(r^{1-\delta_1'}t^{-1/2})^{-N}
 \;\;\;.
 \]
 Choosing $N$ large ensures that this expression is uniformly bounded
 for all $r\geq1$ and $t\leq 1$.
 
 \smallskip
 
 The proof of the theorem is now complete.\hfill$\Box$

 \bigskip
 
 Theorems~\ref{tcsg6.1}  and \ref{tcsg6.2} have a number of standard implications.
 First one may convert the volume bounds of the latter theorem into  Gaussian bounds.
 
 \begin{cor}\label{ccsg6.1}
For each $\varepsilon>0$ there is an $a>0$ such that the  semigroup kernel $K$ of the Gru\v{s}in operator satisfies
 \begin{equation}
0\leq  K_t(x\,;y)\leq a\,(|B(x\,;t^{1/2})|\,|B(y\,;t^{1/2})|)^{-1/2}\,e^{-d(x;y)^2/(4(1+\varepsilon)t)}
\label{ecsg6.21}
\end{equation}
for  all $t>0$ and almost all $x,y\in \Ri^{n+m}$.
\end{cor}

There are several different arguments for passing from on-diagonal kernel bounds to Gaussian bounds
(see, for example, the lecture notes of Grigor'yan \cite{Gri3}).
One proof of the corollary which is in the spirit of the present paper is given by Theorem~4 of 
\cite{Sik3}.
Note that in the latter reference it is implicitly assumed that the kernel is well-defined on the  diagonal
but this is not essential.
One can argue with open sets and near diagonal estimates as in the proof of the comparison theorem,
Theorem~\ref{tcsg3.1}, and in the proof of Theorem~\ref{tcsg6.2}.

\bigskip

It is also a standard argument to pass from the Gaussian bounds of Corollary~\ref{ccsg6.1} and the
conservation property of Theorem~\ref{tcsg6.1} to on-diagonal lower bounds.
But again one has to avoid problems with the definition of the diagonal values.

 \begin{cor}\label{ccsg6.2}
There is an $a>0$ such that the  semigroup kernel $K$ of the Gru\v{s}in operator satisfies
 \begin{equation}
\inf_{X\ni x}|X|^{-2}\int_Xdx\int_Xdy\, K_t(x\,;y)
 \geq a\,|B(x\,;t^{1/2})|^{-1}
 \label{ecsg6.22}
\end{equation}
for  all $t>0$ and almost all $x\in \Ri^{n+m}$ where the average is over open subsets $X$.
\end{cor}
\proof\
At the risk of confusion with the earlier definition we set 
\[
K_t(X\,;Y)=\int_Xdx\int_Ydy\, K_t(x\,;y)=(\one_X,S_t\one_Y)
\]
for each pair of bounded open sets $X,Y$.
Then using self-adjointness and the semigroup property one again verifies that 
\[
|K_t(X\,;Y)|^2\leq K_t(X\,;X)\,K_t(Y\,;Y)
\;\;\;.
\]
Now fix $x$ and $X\ni x$ and  let $Y=B(x\,;R\,t^{1/2})$ for  $R>0$.
Then 
\[
K_t(X\,;Y)=\int_Xdx\Big(1-\int_{Y^{\rm c}}dy\, K_t(x\,;y)\Big)
\]
for all $t>0$ by Theorem~\ref{tcsg6.1}.
But now using the Gaussian bounds of Theorem~\ref{tcsg6.2} and choosing $R$ large one can ensure that
\[
K_t(X\,;Y)\geq |X|/2
\;\;\;.
\]
Then, however, with this choice of $Y$ one has 
\[
|X|^{-2}K_t(X\,;X)\geq (4K_t(Y\,;Y))^{-1}
\;\;\;.
\]
But using the bounds of Theorem~\ref{tcsg6.2} one immediately finds 
\[
K_t(Y\,;Y)\leq a\,\Big(\int_{d(x;y)<Rt^{1/2}}dy\,|B(y\,;t^{1/2})|^{-1/2}\Big)^2 
\;\;\;.
\]
If $y\in B(x\,;R\,t^{1/2})$ the doubling property gives
 \[
 |B(x\,;t^{1/2})|\leq  |B(y\,;d(x\,;y)+t^{1/2})|\leq  |B(y\,;(1+R)t^{1/2})|
 \leq a\,(1+R)^{\widetilde D} |B(y\,;t^{1/2})|
 \]
with $\widetilde D=D\vee D'$.
 In addition 
 \[
 |B(x\,;R\,t^{1/2})|\leq a\,R^{\widetilde D}  |B(x\,;t^{1/2})|
 \;\;\;.
 \]
 Therefore
 \begin{eqnarray*}
K_t(Y\,;Y)&\leq &a\,(1+R)^{2{\widetilde D} }\Big(\int_{d(x;y)<Rt^{1/2}}dy\,|B(x\,;t^{1/2})|^{-1/2}\Big)^2 \\[5pt]
&\leq& a\,(1+R)^{2{\widetilde D} }|B(x\,;t^{1/2})|^{-1} |B(x\,;Rt^{1/2})|^2
\leq a'\,((1+R)R)^{2{\widetilde D} }|B(x\,;t^{1/2})|
\;\;\;.
\end{eqnarray*}
 Combining these estimates gives 
 \[
 |X|^{-2}\int_Xdx\int_Xdy\, K_t(x\,;y)
 \geq a\,|B(x\,;t^{1/2})|^{-1}
 \]
 for all bounded open sets $X$ containing $x$ and this gives the statement of the corollary.
 \hfill$\Box$

\bigskip
\begin{remarkn} \label{rcsg6.1}
If the kernel is continuous then the on-diagonal values are well defined and the corollary gives
 \[
K_t(x\,;x)=\inf_{X\ni x}|X|^{-2}\int_Xdx\int_Xdy\, K_t(x\,;y)
 \geq a\,|B(x\,;t^{1/2})|^{-1}
\]
for all $x\in \Ri^{n+m}$ and $t>0$.
\end{remarkn}

The semigroup kernel $K_t$ of the Gru\v{s}in operator $H$ is not necessarily continuous.
In particular if  $n=1$ and $\delta_1\in[1/2,1\rangle$ then the kernel is  discontinuous. 
The discontinuity is a direct consequence of the fact established in 
\cite{ERSZ1} that the action of the corresponding semigroup $S_t$ on $L_2(\Ri\times \Ri^m)$ is not ergodic.
If $H_+=\{(x_1,x_2)\in\Ri\times \Ri^m: x_1>0\}$ and $H_-=\{(x_1,x_2)\in\Ri\times \Ri^m: x_1<0\}$ then
$S_tL_2(H_+)\subseteq L_2(H_+)$ and $S_tL_2(H_-)\subseteq L_2(H_-)$ for all $t>0$.
This is established as follows.

Let $\varphi\in D(h)\subseteq W^{1,2}(\Ri)$ and  set $\varphi_n=\chi_n\varphi$ with  $\chi_n \colon \Ri \to [0,1]$  defined by
\begin{equation}
\chi_n(x)
= \left\{ \begin{array}{ll}
  0 & \mbox{if } x\leq -1 \;\;\; ,  \\
  {-\log |x|}/{\log n} & \mbox{if } x\in\langle-1, -n^{-1}\rangle \;\;\; , \\
  1 & \mbox{if } x \geq -n^{-1} \;\;\;.\label{ecdeg2.137}
         \end{array} \right.
\end{equation}
Then one verifies that $\|\varphi_n-\one_+\varphi\|_2\to 0$ as $n\to\infty$, where $\one_+$ is the indicator function
of $\overline H_+$, and
$h(\varphi_n-\varphi_m)\to0$ as $n,m\to\infty$.
Thus $\one_+\varphi\in D(\overline h)$.
(See \cite{ERSZ1}, Proposition~6.5 and the discussion in  Section~4 of \cite{RSi}.
Note that it is crucial for the last limit that $\delta\in[1/2,1\rangle$: the conclusion is not valid for 
$\delta\in[0,1/2\rangle$.)
Similarly, replacing $\chi_n(x)$ by $\chi_n(-x)$ one can conclude that 
$\one_{-}\varphi\in D(\overline h)$ with $\one_-$  the indicator function
of $\overline H_-$.
Moreover, $\overline h(\varphi)=\overline h(\one_+\varphi)+\overline h(\one_-\varphi)$.
This suffices to deduce that the semigroup leaves the subspaces $L_2(H_\pm)$ invariant
(see \cite{FOT}, Theorem~1.6.1).

Now one can deduce by contradiction that the kernel has a discontinuity. 
Suppose the kernel $K_t$ of $S_t$ is continuous.
Then it follows from Remark~\ref{rcsg6.1}
that $K_t(x\,;x)\geq a_t>0$ where $a_t=a\,(\sup_{x\in \Ri\times\Ri^m}|B(x\,;t^{1/2}|)^{-1}$.
But since $(\varphi,S_t\psi)=0$ for $\varphi\in L_2(H_+)$ and $\psi\in L_2(H_-)$ one must have
$K_t(x\,;y)=0$ for all $x\in H_+$ and $y\in H_-$.
But this contradicts the continuity hypothesis.

The separation phenomenon  raises the question of boundary conditions  on the hypersurface of separation $x_1=0$.
The closed form $\overline h$ has the decomposition $\overline h(\varphi)=\overline h(\one_+\varphi)+\overline h(\one_-\varphi)$ for all $\varphi\in D(h)$, and by 
closure for all $\varphi\in D(\overline h)$
Since this is the direct analogue  of the decomposition of the form corresponding to the Laplacian
with Neumann boundary conditions on the hypersurface $x_1=0$ it is tempting to describe the separation in terms of Neumann boundary 
conditions.
But this decomposition is misleading since Neumann and  Dirichlet boundary conditions
coincide in this case.
This can be established by potential theoretic reasoning \cite{RSi}.

  First introduce the  form $h_{\rm Dir}$ by restriction of $h$ to the subspace 
  $D(h)\cap L_{2,c}(H_+\oplus H_-)$ where $L_{2,c}(\Omega)$ is defined as the subspace of $L_2(\Omega)$ spanned
  by the functions with compact support.
  Then  $h_{\rm Dir}$ is closable and its closure $\overline h_{\rm Dir}$ corresponds to the operator $H$ with Dirichlet
  boundary conditions imposed at the boundary $x_1=0$.
 Moreover,  $\overline h_{\rm Dir}\geq \overline h$ in the sense of the ordering of forms.
  But comparison of $\overline h_{\rm Dir}$ and $\overline h$ gives a sharp distinction between the
  weakly   degenerate case $\delta_1\in[0,1/2\rangle$ and the strongly degenerate case $\delta_1\in[1/2,1\rangle$.

 \begin{prop}\label{psep}
 Consider the Gru\v{s}in operator with $n=1$.
 If $\delta_1\in[0,1/2\rangle$ then  $\overline h_{\rm Dir}> \overline h$ but if  $\delta_1\in[1/2,1\rangle$ then 
 $\overline h_{\rm Dir}=\overline h$. 
\end{prop}
  \proof\
  First, observe that $D(h)\cap L_{2,c}(H_+\oplus H_-)$ is a core of 
  $\overline h_{\rm Dir}$, by definition, and $\overline h_{\rm Dir}=\overline h$ in restriction to the core.
Secondly, it follows from \cite{RSi}, Proposition~3.2, that  $D(h)\cap L_{2,c}(H_+\oplus H_-)$ is a core of 
  $\overline h$ if, and only if, $C_{\overline h}(\{x_1=0\})=0$ where $C_{\overline h}(A)$ denotes the capacity, with
  respect to $\overline h$ of the measurable set $A$.
  Thus $\overline h_{\rm Dir}=\overline h$ if, and only if, $C_{\overline h}(\{x_1=0\})=0$.
  
Now suppose $\delta\in[1/2,1\rangle$ and define
$\xi_n$  by $\xi_n(x_1)=\chi_n(x_1)\wedge \chi_n(-x_1)$
where $\chi_n$ is given by (\ref{ecdeg2.137}).
Then one verifies that $\xi_n\in [0,1]$ and  $\xi_n(x)=1$  for $x\in[-n^{-1},n^{-1}]$.
Moreover, if $\varphi\in D(h)$ then $\varphi_n=\xi_n\varphi\in D(h)$ and 
 $h(\varphi_n)+\|\varphi_n\|_2^2\to0$
as $n\to\infty$.
But this means that $C_{\overline h}(A)=0$  for each bounded measurable subset of 
the hypersurface $\{x_1=0\}$.
Then it follows by  the monotonicity and  additivity properties of the capacity (see, for example, 
 \cite{FOT}, Section~2.1, or  \cite{BH}, Section~1.8) that  $C_{\overline h}(\{x_1=0\})=0$.
 This establishes the second statement of the proposition.
 
 Finally suppose $\delta\in[0,1/2\rangle$.
 Then one can find a Fourier multiplier $F$ such that $H\geq F$ and, by Proposition~\ref{pcsg3.1},
 one may choose $F$ such that $1+F(p_1,p_2)\geq1+ a\,|p_1|^{2(1-\delta_1)}$ for some $a>0$
 and all $(p_1,p_2)\in \Ri\times \Ri^m$.
 Then if $U\subset \Ri\times\Ri^m$ is an open set with $|U|<\infty$ it follows by the calculation at the end
 of Section~3 of \cite{RSi} that 
 \[
 C_{\overline h}(U)\geq |U|^2(\one_U,(I+F)^{-1}\one_U)^{-1}
 \;\;\;.
 \]
 Now set $U_\varepsilon=\langle-\varepsilon,\varepsilon\rangle\times V$ where $V$ is an open subset of $\Ri^m$.
 Then calculating as in the proof of Proposition~4.1 of \cite{RSi} one finds
 \begin{eqnarray*}
  C_{\overline h}(U_\varepsilon)&\geq& 4\,\varepsilon^2\,|V|^2\,\bigg(4\,\varepsilon^2\int_{\Ri^m}dp_2\,(\tilde\one_V(p_2))^2
  \int_\Ri dp_1\,(1+F(p_1,p_2))^{-1}(\sin(\varepsilon p_1)/(\varepsilon p_1))^2\bigg)^{-1}\\[5pt]
  &\geq& |V|\,\bigg(\int_\Ri dp_1\,\left(1+a\,|p_1|^{2(1-\delta_1)}\right)^{-1}\bigg)^{-1}\geq a_{\delta_1}\,|V|
  \end{eqnarray*}
 where $a_{\delta_1}>0$.
 Note that the strict positivity of $a_{\delta_1}$ requires $\delta_1\in[0,1/2\rangle$.
Then,  however, one must have  $C_{\overline h}(\{x_1=0\})>0$.\hfill$\Box$

\bigskip

The moral of the proposition  is that 
in the strongly degenerate case $\delta\in[1/2,1\rangle$ the Dirichlet and Neumann boundary conditions coincide.
The separation is a spontaneous effect which is not characterized by a particular choice of boundary conditions.



 \section{A one dimensional example}\label{Scsg7}
 
In this section we give a further analysis of the one-dimensional example
 discussed in \cite{ERSZ1}, Sections~5 and 6.
 This example is a special case of the Gru\v{s}in operator with $n=1$,  $m=0$ and $\delta_1'=0$.
 Its structure provides a guide to the anticipated structure of the more interesting examples with
  $n=1$ and $m\geq1$.

Let  $h(\varphi)=(\varphi',c_\delta\,\varphi')$ be the form on $L_2(\Ri)$ 
 with domain $W^{1,2}(\Ri)$ where $c_\delta$ is given by $c_\delta(x)=(x^2/(1+x^2))^\delta$ with $\delta>0$.
 The form is closable  and we let $\overline h$ denote its  closure (relaxation).
 Let $H$ be  the corresponding positive self-adjoint
 operator on $L_2(\Ri)$,
  $S$ the semigroup generated by $H$  and $K$ the   kernel of $S$.
 All the qualitative features we subsequently derive extend to the semigroups associated with forms
 $h(\varphi)=(\varphi',c\,\varphi')$ with $c\in L_\infty(\Ri)$ and $c\sim c_\delta$.
 
The  Riemannian distance is now  given by
$d(x\,;y)=|\int^x_ydt\,c_\delta(t)^{-1/2}|$.
If $\delta\in[0,1\rangle$ then $d(\cdot\,;\cdot)$ is a genuine distance but if $\delta\geq 1$
 then the distance between the left and half right lines is infinite.
We concentrate on the  case $\delta\in[0,1\rangle$   and comment on the distinctive features of the case $\delta\geq 1$
at the end of the section.

If  $\delta\in[0,1\rangle$ the Riemannian ball $B(x\,;r)=\{y:d(x\,;y)<r\}$ is an interval
and  the volume is the length of the interval.
It is straightforward to estimate this length from the explicit form of $c_\delta$.
One finds
$|B(x\,;r)|\sim  r $  if  $|x|\geq 1$ or if $r\geq 1$,
$ |B(x\,;r)|\sim r^{1/(1-\delta)}$  if $|x|\leq 1$, $r\leq 1$ and $d(0\,;x)<r$ and 
$ |B(x\,;r)|\sim |x|^\delta\, r$
if $|x|\leq 1$, $r\leq 1$ and $d(0\,;x)\geq r$.
Then $|B|$ satisfies the doubling property (\ref{epre1.17})
with doubling dimension  ${\widetilde D} =D=1/(1-\delta)$.
These bounds are all consistent with the general estimates of Proposition \ref{pcsg5.1}.

Moreover, if $B_+(x\,;r)=B(x\,;r)\cap \Ri_+$   then $|B(x\,;r)|/2\leq |B_+(x\,;r)|\leq |B(x\,;r)|$
  for $x\geq0$. 
  Hence $|B_+|$ satisfies similar estimates for $x\geq0$.

\medskip

Now we consider bounds on the associated semigroup kernel $K_t$.
There are two distinct cases $\delta\in[0,1/2\rangle$ and $\delta\in[1/2,1\rangle$.
These correspond to weak and strong degeneracy in the terminology of \cite{ABCF} and \cite{MaV}
who have considered control theory aspects of similar one-dimensional examples.

\smallskip

\noindent{\bf Case I} $\;\delta\in[0,1/2\rangle$.
The restriction on $\delta$ means that the degeneracy of the coefficient is relatively mild and  the general conclusions of \cite{Tru2} \cite{FKS} are applicable.
But in the present situation one can deduce much more.
The Gaussian upper bounds of Corollary~\ref{ccsg6.1}
 are valid but there are matching lower bounds.

\begin{prop}\label{pcsg4.1}
If $\delta\in[0,1/2\rangle$  there are $b,c>0$ such that 
\begin{equation}
K_t(x\,;y)\geq b\,|B(x\,;t^{1/2})|^{-1}e^{-cd(x;y)^2/t}
\label{ep2}
\end{equation}
for all $x,y\in\Ri$ and $t>0$.
\end{prop}

The deduction of the  off-diagonal lower bounds requires  some additional information
on continuity.

\begin{lemma}\label{lcsg4.4}
If $\delta\in[0,1/2\rangle$ then there is an $a>0$ such that
\[
|\varphi(x)-\varphi(y)|^2\leq a\,{{d(x\,;y)^2}\over{V(x\,;y)}}\,h(\varphi)
\]
for all $x,y\in\Ri$ and all $\varphi\in D(h)$ where $V(x\,;y)=|B(x\,;d(x\,;y))|\vee |B(y\,;d(x\,;y))|$.
\end{lemma}
\proof\
Let $\varphi\in W^{1,2}(\Ri)$.
Then
\begin{eqnarray*}
|\varphi(x)-\varphi(y)|^2&=&\Big|\int^y_xds\,\varphi'(s)\Big|^2\\[5pt]
   &\leq&\Big|\int^y_xds\,c_\delta(s)^{-1}\Big|\,\Big|\int^y_xds\,c_\delta(s)\,\varphi'(s)^2\Big|
     \leq\Big|\int^y_xds\,c_\delta(s)^{-1}\Big|\,h(\varphi)
     \;\;\;.
     \end{eqnarray*}
Now the proof of the lemma follows from the upper bound in the next lemma.

\smallskip

Note that at this point it is essential that $\delta\in[0,1/2\rangle$ to ensure that $c_\delta^{-1}$ 
is locally integrable.

\begin{lemma}\label{lcsg4.5}
If $\delta\in[0,1/2\rangle$ then there is an $a>0$ such that
  \[
{{d(x\,;y)^2}\over{|x-y|}}\leq \Big|\int^y_xds\,c_\delta(s)^{-1}\Big|
\leq a\,{{d(x\,;y)^2}\over{V(x\,;y)}}
\]
for all $x,y\in\Ri$.
\end{lemma}
\proof\
The left hand bound follows directly from the Cauchy--Schwarz inequality;
\[
d(x\,;y)^2=\Big|\int^y_xds\,c_\delta(s)^{-1/2}\Big|^2
\leq |x-y|\,\Big|\int^y_xds\,c_\delta(s)^{-1}\Big|
\;\;\;.
\]
The right hand bound uses the volume estimates and follows by treating various different cases.

First  set $D(x\,;y)=|\int^y_xds\,c_\delta(s)^{-1}|$.
Then if $|x-y|\geq 1/2$ one has $d(x\,;y)\sim |x-y|$ and by similar reasoning
 $D(x\,;y)\sim |x-y|$.
Moreover, $V(x\,;y)\sim|x-y|$ by  the volume estimate with $r\geq 1$.
Therefore  $D(x\,;y)\sim d(x\,;y)^2/V(x\,;y)$.
If, however, $x,y\geq 1/2$ or $x,y\leq-1/2$ then the estimate follows by similar reasoning
but using the volume estimate with $|x|\geq 1$.
It remains to consider  $x,y$ such that $|x-y|\leq 1$.
But  by symmetry the discussion can be reduced to
 two cases $0\leq y<x\leq 1$ and $-1\leq y<0<x\leq 1$.

 Consider the first case. 
 Then $D(x\,;y)\sim x^{1-2\delta}-y^{1-2\delta}$ and 
$d(x\,;y)\sim x^{1-\delta}-y^{1-\delta}$ by explicit calculation.
 Moreover,  $V(x\,;y)=|B(x\,;d(x\,;y))|\leq a\,(x-y)$ for some $a\geq1$.
 But
\begin{eqnarray*}
(x^{1-\delta}-y^{1-\delta})^2-(x^{1-2\delta}-y^{1-2\delta})(x-y)
=(x^{1/2-\delta}y^{1/2}-x^{1/2}y^{1/2-\delta})^2\geq0
\;\;\;.
\end{eqnarray*}
This again establishes the required bound.

Finally consider the second case and suppose that $x\geq |y|$. 
 Then $d(x\,;y)\sim d(x\,;0)\sim x^{1-\delta}$ and $D(x\,;y)\sim D(x\,;0)\sim x^{1-2\delta}$.
 Moreover,
 \[
 V(x\,;y)=|B(x\,;d(x\,;y))|\leq |B(x\,;2d(x\,;0))|\leq a\,|B(x\,;d(x\,;0))|\leq a'\,x
 \]
 where the second estimate uses volume doubling.
 The required bound follows immediately.
 The case $|y|>x$ is similar with the roles of $x$ and $y$ interchanged.
 \hfill$\Box$
 
 \bigskip

The  bound in Lemma~\ref{lcsg4.4} is now an immediate consequence of the upper bound of Lemma~\ref{lcsg4.5}.
Moreover the Gaussian lower bound in Proposition~\ref{pcsg4.1}
 follows directly from the  the Gaussian upper bound and the continuity bound of Lemma~\ref{lcsg4.4}. This last implication is, for example, a direct consequence of Theorem~3.1 in \cite{Cou4}.
 This theorem is applied with $w=2$, $p=2$ and  $\alpha=1$. 
 Note that the doubling dimension
$D=1/(1-\delta)< \alpha \,p=2$ because $\delta\in[0,1/2\rangle$.
Therefore Theorem~3.1 of \cite{Cou4} is indeed applicable.
This completes the proof of Proposition~\ref{pcsg4.1}.\hfill$\Box$
 
\bigskip

\noindent{\bf Case 2} $\;\delta\in[1/2,1\rangle$.
This case corresponds to strong degeneracy in the terms of \cite{ABCF} and \cite{MaV} and 
describes a quite different situation. 
It follows from Proposition~6.5 of \cite{ERSZ1} that $S_tL_2(\Ri_\pm)\subseteq L_2(\Ri_\pm)$.
Thus the system separates into two ergodic components and the semigroup kernel has the property $K_t(x\,;y)=0$
for $x<0$ and $y>0$.
One can, however, extend the foregoing analysis to the two components.
First we prove that the  kernel of the semigroup restricted to $L_2(\Ri_+)$
is  H\"older continuous.
Note that the generator of the restriction of the semigroup to $L_2(\Ri_+)$ is the operator
associated with the  closure of the form  obtained by restricting $\overline h$  to 
$W^{1,2}(\Ri_+)=\one_{\Ri_+}W^{1,2}(\Ri)$.
Now the continuity  proof is by a variation of the usual   Sobolev inequalities
\[
|\varphi(x)|^2 \leq a\,(\|\varphi'\|_2^2 +\|\varphi\|_2^2)
\]
and
\[
|x-y|^{-2\gamma}|\varphi(x)-\varphi(y)|^2 \leq a_\gamma \,(\|\varphi'\|_2^2 +\|\varphi\|_2^2)
\]
where  $\gamma\in\langle0,1/2\rangle $.
Fix $\sigma>0$ and let  $\psi \in C_c^\infty(\sigma /2,\infty)$ with 
$\psi(x)=1$ for $x>\sigma $.
Then
\begin{eqnarray*}
\|(\psi \varphi)'\|_2^2 +\|\psi \varphi\|_2^2&\leq&
\|\psi' \varphi+\varphi'\psi\|_2^2+\|\psi \varphi \|_2^2 \leq 2\,\|\psi' \varphi\|_2^2+
2\,\|\psi\varphi'\|_2^2 +\|\psi\varphi\|_2^2\\[5pt]
&\leq& 2\,\|\psi\varphi'\|_2^2+(1+4/\sigma )\,\| \varphi\|_2^2\leq 2\,
c_\delta(\sigma /2)^{-1}\|c_\delta^{1/2}\varphi'\|_2^2+(1+4/\sigma )\,\| \varphi\|_2^2\\[5pt]
&\le& a_\sigma \,(h(\varphi)+\|\varphi\|_2^2)\;\;\;.
\end{eqnarray*}
Therefore
\[
|\varphi(x)|^2=|(\psi\varphi)(x)|^2 \leq a'_\sigma \,(h(\varphi)+\|\varphi\|_2^2)
\]
and
\[
|x-y|^{-2\gamma}|\varphi(x)-\varphi(y)|^2 \leq a_{\sigma ,\gamma} \,(h(\varphi)+\|\varphi\|_2^2)
\]
for all $x, y\geq \sigma $ and all $\varphi\in W^{1,2}(\Ri_+)$.
These bounds then extend by continuity to the closure of  $h$ and are sufficient to deduce
that the semigroup kernel $K_t$ is uniformly bounded and H\"older continuous on 
$[\sigma ,\infty\rangle\times [\sigma ,\infty\rangle$ for each $\sigma >0$.
In particular $K_t$ is continuous on $\langle 0,\infty\rangle\times \langle 0,\infty\rangle$.
But the bounds depend on $\sigma $ and  do not give good {\it a priori} bounds on 
$\sup_{x,y>0}K_t(x\,;y)=\|S_t\|_{1\to\infty}$.
This can again be accomplished by a Nash inequality argument.

If $h_+$ temporarily denotes the closed form of the generator of the semigroup on $L_2(\Ri_+)$ and 
$E\varphi$ denotes the symmetric extension of $\varphi\in L_2(\Ri_+)$ to $E\varphi\in L_2(\Ri)$ then 
${\overline h}(E\varphi)=2\,h_+(\varphi)$ for all $\varphi\in D(h_+)$.
Moreover, $\|E\varphi\|_2^2=2\,\|\varphi\|_2^2$ and $\|E\varphi\|_1^2=4\,\|\varphi\|_1^2$.
Therefore the Nash inequalities (\ref{epre1.7}) give similar inequalities
\[
\|\varphi\|_2^2\leq r^{-2}h_+(\varphi)
+\pi^{-1}\, V_F(r)\,\|\varphi\|_1^2
\]
for all $\varphi\in D(h_+)$.
Hence the  semigroup restricted to $L_2(\Ri_+)$, or  $L_2(\Ri_-)$, again satisfies bounds
 $\|S_t\|_{1\to\infty}\leq a\,t^{(-1/(2(1-\delta)),-1/2)}$ for all $t>0$.
In particular, the kernel $K_t$ is uniformly bounded on $\Ri_+\times \Ri_+$,
or  $\Ri_-\times \Ri_-$,
for each $t>0$.

One can then prove the analogues of Corollaries~\ref{ccsg6.1}
and ~\ref{ccsg6.2} on the half-lines.
\begin{prop}\label{pcsg4.2}
If $\delta\in[1/2,1\rangle$ and  $B_\pm(x\,;r)=B(x\,;r)\cap \Ri_\pm$ then
for each $\varepsilon\in\langle0,1]$ there is an $a>0$ such that 
\begin{equation}
K_t(x\,;y)\leq a\,|B_\pm(x\,;t^{1/2})|^{-1}e^{-d(x;y)^2/(4t(1+\varepsilon))}
\label{ep3}
\end{equation}
for all $x,y\in\Ri_\pm$ and $t>0$.
Moreover, there is  $b>0$ such that 
\begin{equation}
K_t(x\,;x)\geq b\,|B_\pm(x\,;t^{1/2})|^{-1}\label{ep4}
\end{equation}
for all $x\in\Ri_\pm$ and $t>0$.
\end{prop}

In the one-dimensional case the statement of Proposition~\ref{psep}
 can also be described in terms of the vector fields defining $h$ and $h_{\rm Dir}$.

 Let $X=c_\delta^{1/2}\,d$, with $d=d/dx$,  denote the $C^\delta$-vector field acting on $L_2(\Ri)$  with domain $D(X)=C_c^\infty(\Ri)$  and $X_0$ the restriction of $X$ to 
 $D(X_0)=C_c^\infty(\Ri\backslash\{0\})$.
 Then $h(\varphi)=\|X\varphi\|_2^2$ and $h_{\rm Dir}(\varphi)=\|X_0\varphi\|_2^2$.
 It follows straightforwardly that the corresponding positive self-adjoint operators are given by
 $H=X^*\overline X$ and $H_{\rm Dir}=X_0^*\overline X_0$.
 But Proposition~\ref{psep} can now be restated as follows.
 
 \begin{cor}\label{csep}
    If $\delta\in[0,1/2\rangle$ then  $\overline X_0\subset \overline X$ but if  $\delta\in[1/2,1\rangle$ then  $\overline X_0= \overline X$.
    \end{cor}

Although we have restricted attention to the case $\delta\in[0,1\rangle$ the form $h$ is densely defined
and closable for all $\delta\geq 0$.
But the semigroup and its kernel have very different properties if  $\delta\geq 1$.
The properties of the distance
$d(x\,;y)=|\int^x_ydt\,c_\delta(t)^{-1/2}|$ are also quite different.
The distance is finite on the open half line $\langle 0,\infty\rangle$ but since the integral diverges at zero
the distance to the origin is at infinity.
Now consider points $x,y\in\langle0,1]$. Then $d(x\,;y)\sim |\ln x/y\,|$ if $\delta=1$ and $d(x\,;y)\sim|x^{1-\delta}-y^{1-\delta}|$
if $\delta>1$.
Consider the case $\delta=1$ with the equivalent distance $\tilde d(x\,;y)=|\ln x/y\,|$. 
Let $x_n\in\langle0,1\rangle$ be an arbitrary sequence which converges downward to zero.
One may assume $x_1\leq e^{-1}$.
Then $\tilde d(x_n\,;x_ne)=1=\tilde d(x_n\,;x_ne^{-1})$ and $|B(x_n\,;1)|=2\,x_n\, \sinh 1$.
Therefore $|B(x_n\,;1)|\to 0$ as $n\to\infty$.

Alternatively if $x_n=e^{-n}$ with $n\geq 1$ then $B(x_n\,;n/2)=\langle e^{-3n/2},e^{-n/2}\rangle$ and 
$B(x_n\,;n)=\langle e^{-2n},1\rangle$. 
Therefore $|B(x_n\,;n/2)|=e^{-n/2}(1-e^{-n})\to0$ and  
$|B(x_n\,;n)|=1-e^{-2n}\to1$ as $n\to\infty$ so the volume cannot satisfy the volume
doubling property.

These divergences allow one to argue that the semigroup kernel is not bounded near the origin.

\smallskip

 \section{Applications}\label{Scsg8}
 
 In the foregoing we  established that Gru\v{s}in operators have many important properties
 in common with strongly elliptic operators; the wave equation has a finite propagation speed,
 the heat kernel satisfies Gaussian upper bounds and the heat semigroup conserves probability.
 Then one can readily adapt arguments developed for strongly elliptic operators to obtain
 further detailed information about the Gru\v{s}in operators, e.g.,  information on boundedness of Riesz 
 transforms, spectral multipliers and Bochner--Riesz summability, holomorphic functional calculus, Poincar\'e inequalities and maximal   regularity.
 We conclude by describing briefly some of these applications to Gru\v{s}in operators which are 
 a straightforward consequence of general theory and which require no further detailed arguments.
 It should, however, be emphasized that there are significant differences between the degenerate and the non-degenerate  theories  related to continuity and positivity properties.
 In particular the one-dimensional example in Section~\ref{Scsg7} demonstrates that 
 the heat kernel  is not necessarily continuous nor strictly positive.
 Therefore there are limitations to possible extensions of the results of classical analysis to the degenerate case.

\subsection{Boundedness of Riesz transforms}

First we consider boundedness of the Riesz  transforms associated with a general Gru\v{s}in operator $H$ on $L_p(\Ri^{n+m})$ for $p\in\langle1,2] $. 
The result can be stated  in terms of the {\bf carr\'e du champ} associated with $H$
(see, for example, Section~I.4 of \cite{BH}). 
Formally  the carr\'e du champ is given by 
\[
\Gamma_\psi=\psi(H\psi)-2^{-1}H\psi^2
\;\;\;.
\]
Note that if $\psi \in C_c^\infty(\Ri^{n+m})$ then
\[
 \Gamma_\psi(x)=\sum_{i,j=1}^{n+m}c_{ij}(x)(\partial_i \psi)(x)( \partial_j\psi)(x)
\]
and
\[
\| \Gamma_\psi\|_1=\|H^{1/2}\psi\|_2
\;\;\;.
\]
Now one can formulate the
result concerning boundedness of  the   Riesz
transform in an analogous manner.
\begin{thm}\label{tcsgc.1}
If $H$ is   Gru\v{s}in  operator then
$$
\|\Gamma_{H^{-1/2}\psi}\|_{p/2}\leq \|\psi\|_{p}
$$
for all $\psi \in L^p(\Ri^{n+m})$ and all
$p\in\langle1,2]$.
In addition the map $\psi \to
\Gamma_{H^{-1/2}\psi}^{1/2}$ is  weak type~$(1,1)$, i.e.,
\[
|\{ x \in X : |  \Gamma_{H^{-1/2}\psi} (x)|^{1/2} >
\lambda \}
\leq a \, {\| \psi \|_{1}}/{\lambda}
\]
for all $ \lambda \in \Ri_+$ and all  $ \psi \in L_1(\Ri^{n+m})$.
\end{thm}
\proof\
The proof of Theorem~\ref{tcsgc.1} is a straightforward modification of
the proof of Theorem~5 of \cite{Sik3}.
The assumptions of Theorem~5 of \cite{Sik3} hold in virtue of the property of
finite speed of propagation proved
in Proposition~4.1 and the kernel bounds (\ref{ecsg6.2}) of Theorem~\ref {tcsg6.2} .
\hfill$\Box$

\subsection{Spectral multipliers}
Each  Gru\v{s}in operator  $H$ is  positive definite and
self-adjoint.
Therefore $H$ admits a
spectral  resolution  $E_H(\lambda)$  and for  any  bounded  Borel
function $F\colon [0, \infty) \to \Ci$ one can define the operator $F(H)$
by 
   \begin{equation}\label{equw}
   F(H)=\int_0^{\infty}d E_H(\lambda)\,F(\lambda) \;\;\;.
   \end{equation}
It then follows that $F(H)$ is bounded on
$L^2(\Ri^{n+m})$.
Spectral  multiplier  theorems  investigate  sufficient conditions on
function  $F$ which  ensure  that the   operator  $F(H)$ extends  to a
bounded operator on  $L_q$ for some $q\in[1,\infty]$.
\begin{thm}\label{tcsgc.2}
If $H$ is a  Gru\v{s}in  operator,
$s > (D\vee D')/2$ and $F \colon [0,\infty) \to \Ci$
 is  a  bounded Borel function  such that
     \begin{equation}\label{Ale}
     \sup_{t>0}\| \eta \,\delta_t F \|_{W^{s,\infty}} < \infty,
     \end{equation}
   where $ \delta_t F(\lambda)=F(t\lambda)$ and
   $\| F \|_{W^{s,p}}=\|(I-d^2/d x^2)^{s/2}F\|_{L_p}$. 
   Then $F(H)$ is
 weak type $(1,1)$ and bounded on~$L_q$ for all  $q\in\langle1,\infty\rangle$.
\end{thm}
\proof\
The proof of Theorem~\ref{tcsgc.2} is a direct consequence of the Gaussian bounds  (\ref{ecsg6.21}) on the heat kernel corresponding to $H$
given by Corollary~\ref{ccsg6.1} and Theorem~3.1 of \cite{DOS}.
(See also Theorem~3.5 of \cite{CS}.)
\hfill$\Box$

\bigskip

The theory of spectral multipliers is  related to  and motivated
by  the study of convergence  of the Riesz   means or  convergence   of other
eigenfunction expansions of self-adjoint  operators.
To  define  the
Riesz means of the operator $H$ we set
   \begin{equation}\label{vab1}
   \sigma^{s}_R(\lambda)=
      \left\{
       \begin{array}{cl}
       (1-\lambda/R)^{s}  &\mbox{for}\;\; \lambda \le R \\
       0  &\mbox{for}\;\; \lambda > R. \\
       \end{array}
      \right.
   \end{equation}
We then define  $\sigma^{s}_R(H)$  by spectral theory.
The operator $\sigma^{s}_R(H)$ is the Riesz or the Bochner-Riesz mean of
order $s$. 
The basic question in the theory of the Riesz means is  to
establish  the critical exponent for continuity and convergence of
the Riesz means.
More precisely one wishes to ascertain  the optimal range of
$s$ for which the Riesz means $\sigma^{s}_R(H)$ are uniformly
bounded on $L_1(\Ri^{n+m})$. 
A result of this type is given by the
following.
\begin{thm}\label{tcsgc.3}
If $H$ is a  Gru\v{s}in  operator and     $s > (D\vee D')/2$
then
\[
   \sup_{R>0}\|\sigma_R^{s}(H)\|_{q \to q}
     \leq a < \infty
\]
 for all $q\in [1,\infty]$.
   Hence  
\[
     \lim_{R \to \infty}
     \|\sigma_R^{s}(H)\varphi-\varphi\|_{q \to q}=0
\]
for all $q\in [1,\infty]$ and all $\psi\in L_q(\Ri^{n+m})$.
\end{thm}
\proof\
The proof of Theorem \ref{tcsgc.3} is again a  direct consequence of
the Gaussian bounds (\ref{ecsg6.21}) on the heat kernel corresponding to $H$
and
Corollary~6.3 of \cite{DOS}.
\hfill$\Box$

\bigskip

Next we consider the implication of  the Gaussian bounds on the heat kernel
for the holomorphic
function calculus of the Gru\v{s}in operators.
First we briefly recall the notion of holomorphic
function calculus.
For each $\theta >  0$ set  $\Sigma(\theta)=\{z\in C\backslash\{0\}\colon
|\mbox{arg}\,z|<\theta$\}. Let $F$ be a
bounded   holomorphic  function  on     $\Sigma(\theta).$
By   $\Vert F \Vert_{\theta,\infty}$ we denote the supremum of $F$ on
$\Sigma(\theta)$.
The general problem of interest is to  find sharp bounds, in terms of
$\theta$, of the norm of $F(H)$ as an  operator  acting on
$L_p(\Ri^{n+m})$.
  It is known (see  \cite{CDMY}, Theorem~4.10) that
these bounds on  the   holomorphic
functional calculus when $\theta$ tends to~$0$ are related to
spectral multiplier theorems for $H$. 
The following theorem describing
  holomorphic function calculus for  Gru\v{s}in  type operators
follows from (\ref{ecsg6.21}) and Corollary~\ref{ccsg6.1}.
\begin{thm}\label{tcsgc.4}
If  $H$ is a  Gru\v{s}in  operator and 
$s> (D\vee D')|1/p-1/2|$
then
\[
\|  F(H) \|_{p \to p} \le {a}\,{\theta^{-s}}
\| F  \|_{\theta, \infty}
\]
for all $\theta > 0$.
\end{thm}
\proof\
Theorem \ref{tcsgc.4} follows from
the Gaussian bounds (\ref{ecsg6.21}) of
Corollary~\ref{ccsg6.1} and
Proposition~8.1 of \cite{DOS}.
\hfill$\Box$

\subsection{Concluding remarks and comments}

The above statements on the boundedness of the Riesz transforms
and the spectral multipliers for Gru\v{s}in  operators  are not always optimal. 
Using the basic estimates of 
Corollary~\ref{ccsg6.1} and Proposition~4.1 one can analyze the
boundedness
of the Riesz transforms for $p>2$.

In the multiplier result  discussed above, Theorem~\ref{tcsgc.2}, the
critical
exponent required for the  order of differentiability of the function
$F$ is equal to  half of the homogeneous dimension $D\vee D'$.
This  is a  quite typical situation and  for the standard Laplace operator
this exponent is optimal. 
We expect,  however,  that in many cases  it is
possible to obtain multiplier results for Gru\v{s}in  operators with
critical exponent essentially smaller then the half of the homogeneous
dimension $D\vee D'$.

It is also possible to obtain a  version of the Poincare inequality and Nash type
results
similar to those discussed in Section~\ref{Scsg7}. 
But  results of this nature require
substantial new proofs which we hope to describe elsewhere
We conclude by  stressing that
Corollary~\ref{ccsg6.1} and Proposition~\ref{pscsg4.0}  provide a sound basis for
further analysis of Gru\v{s}in
type operators.

\subsection*{Acknowledgement}
This work  was  supported by an Australian Research Council (ARC) Discovery Grant DP 0451016.
It grew out of an earlier collaboration with Tom ter Elst to whom the authors are indebted for many
helpful discussions about degenerate operators.

\end{document}